\documentclass[11pt,a4paper]{article}
\usepackage{amsfonts,amsmath,amssymb,hyperref,accents,amsthm}
\usepackage{tikz}
\usetikzlibrary{calc,arrows,positioning,fit,petri}

\newcommand{\dis}{\displaystyle}
\newcommand{\txt}{\textstyle}

\newcommand{\noi}{\noindent}
\newcommand{\halmos}{\rule{1ex}{1.4ex}}
\newcommand{\QED}{\nopagebreak{\hspace*{\fill}$\halmos$\medskip}}

\newcommand{\med}{\medskip}
\newcommand{\quand}{\quad\mbox{and}\quad}

\newtheoremstyle{mythm}
  {}
  {}
  {\itshape}
  {}
  {\bfseries}
  {}
  {.5em}
  {#1 #2 \thmnote{(#3)}}

\theoremstyle{mythm}
\newtheorem{theorem}{Theorem}
\newtheorem{proposition}[theorem]{Proposition}
\newtheorem{lemma}[theorem]{Lemma}
\newtheorem{exercise}[theorem]{Exercise}
\newtheorem{corollary}[theorem]{Corollary}
\newtheorem{conjecture}[theorem]{Conjecture}

\newtheorem{counterex}[theorem]{Counterexample}

\newcommand{\bt}{\begin{theorem}}
\newcommand{\et}{\end{theorem}}
\newcommand{\bl}{\begin{lemma}}
\newcommand{\el}{\end{lemma}}
\newcommand{\bp}{\begin{proposition}}
\newcommand{\ep}{\end{proposition}}
\newcommand{\bcor}{\begin{corollary}}
\newcommand{\ecor}{\end{corollary}}
\newcommand{\br}{\begin{remark}\rm}
\newcommand{\er}{\end{remark}}
\newcommand{\bcon}{\begin{conjecture}}
\newcommand{\econ}{\end{conjecture}}
\newcommand{\bex}{\begin{exercise}}
\newcommand{\eex}{\end{exercise}}
\newcommand{\bcou}{\begin{counterex}}
\newcommand{\ecou}{\end{counterex}}

\newenvironment{Proof}[1][]{\noi\textbf{Proof #1}}{\QED}
\newcommand{\bpro}{\begin{Proof}}
\newcommand{\epro}{\end{Proof}}

\newcommand{\be}{\begin{equation}}
\newcommand{\ee}{\end{equation}}
\newcommand{\ba}{\begin{array}}
\newcommand{\ea}{\end{array}}
\newcommand{\bc}{\be\begin{array}{r@{\,}c@{\,}l}}
\newcommand{\ec}{\end{array}\ee}
\newcommand{\bac}{\begin{array}{r@{\,}c@{\,}l}}

\newcommand{\al}{\alpha}

\newcommand{\ga}{\gamma}
\newcommand{\Ga}{\Gamma}
\newcommand{\de}{\delta}

\newcommand{\eps}{\varepsilon}
\newcommand{\la}{\lambda}
\newcommand{\La}{\Lambda}
\newcommand{\sig}{\sigma}

\newcommand{\si}{\ensuremath{\sigma}}

\newcommand{\Di}{{\cal D}}
\newcommand{\Ei}{{\cal E}}

\newcommand{\Gi}{{\cal G}}

\newcommand{\Mi}{{\cal M}}

\newcommand{\Ri}{{\cal R}}
\newcommand{\Si}{{\cal S}}

\newcommand{\Wi}{{\cal W}}
\newcommand{\Xc}{{\cal X}}

\newcommand{\R}{{\mathbb R}}

\newcommand{\Z}{{\mathbb Z}}

\renewcommand{\P}{{\mathbb P}}
\newcommand{\E}{{\mathbb E}}

\newcommand{\desd}{\ensuremath{\Leftrightarrow}}
\newcommand{\volgt}{\ensuremath{\Rightarrow}}
\newcommand{\up}{\uparrow}
\newcommand{\down}{\downarrow}
\newcommand{\sub}{\subset}
\newcommand{\beh}{\backslash}

\newcommand{\asto}[1]{\underset{{#1}\to\infty}{\longrightarrow}}

\newcommand{\ti}{\tilde}

\newcommand{\ov}{\overline}

\newcommand{\ffrac}[2]{{\textstyle\frac{{#1}}{{#2}}}}

\newcommand{\di}{\mathrm{d}}
\newcommand{\half}{{[0,\infty)}}
\newcommand{\expo}{\mbox{\large\it e}}
\newcommand{\ex}[1]{\expo^{\,\textstyle{#1}}}

\newcommand{\ha}{\ffrac{1}{2}}

\newcommand{\Ly}{V}
\newcommand{\Vv}{u}
\newcommand{\cri}{{\rm c}}

\begin{document}

\makeatletter\@addtoreset{equation}{section}
\makeatother\def\theequation{\thesection.\arabic{equation}}

\renewcommand{\labelenumi}{{\rm (\roman{enumi})}}
\renewcommand{\theenumi}{\roman{enumi}}

\title{Rigorous results for the Stigler-Luckock model\\
for the evolution of an order book}
\author{Jan M. Swart\vspace{6pt}\\
{\small Institute of Information Theory and Automation of the ASCR (UTIA)}\\
{\small Pod vod\'arenskou v\v{e}\v{z}\'i 4}\\
{\small 18208 Praha 8}\\
{\small Czech Republic}\\
{\small e-mail: swart@utia.cas.cz}
\vspace{4pt}}
\date{\today}
\maketitle

\begin{abstract}\noi
In 1964, G.J.~Stigler introduced a stochastic model for the evolution of an
order book on a stock market. This model was independently rediscovered and
generalized by H.~Luckock in 2003. In his formulation, traders place buy and
sell limit orders of unit size according to independent Poisson processes with
possibly different intensities. Newly arriving buy (sell) orders are either
immediately matched to the best available matching sell (buy) order or stay in
the order book until a matching order arrives. Assuming stationarity, Luckock
showed that the distribution functions of the best buy and sell order in the
order book solve a differential equation, from which he was able to calculate
the position of two prices $J^\cri_-<J^\cri_+$ such that buy orders below
$J^\cri_-$ and sell orders above $J^\cri_+$ stay in the order book forever
while all other orders are eventually matched. We extend Luckock's model by
adding market orders, i.e., with a certain rate traders arrive at the market
that take the best available buy or sell offer in the order book, if there is
one, and do nothing otherwise. We give necessary and sufficient conditions for
such an extended model to be positive recurrent and show how these conditions
are related to the prices $J^\cri_-$ and $J^\cri_+$ of Luckock.
\end{abstract}

\vspace{.4cm}
\noi
%
{\it Acknowledgments.} Work sponsored by GA\v{C}R grants 16-15238S and
P201/12/2613.


{\setlength{\parskip}{-2pt}\tableofcontents}

\section{Introduction and results}

\subsection{Definition of the model}\label{S:intro}

We will be interested in a stochastic model for traders interacting through an
order book as is commonly used on a stock market or commodity market. In the
more theoretical economic literature, the sort of trading system we are
interested in is also known as the continuous double auction. In our specific
model of interest, traders arrive according to independent Poisson processes
and place either a buy or sell limit order for exactly one item of a certain
stock or commodity. If the order book already contains a suitable offer, then
the new limit order is immediately matched with the best available offer,
i.e., a new buy limit order at a price $x$ is cancelled against an existing
sell limit order at the lowest possible price $x'\leq x$, if such a sell limit
order exists, and vice versa for new sell limit orders. Orders that are not
immediately matched stay in the order book until they are matched with a new
incoming order, or, if such an order never comes, forever. This model, in
discrete time and for a specific choice of the parameters, was invented by
Stigler \cite{Sti64} and, in its full generality, independently by Luckock
\cite{Luc03}. The model was subsequently again independently reinvented by
Pla\v{c}kov\'a in her master thesis \cite{Pla11} and by Yudovina in het
Ph.D.\ thesis \cite{Yud12a,Yud12b}. We will generalize the
model by also allowing market orders, i.e., with a certain rate a trader
arrives that takes the best available limit buy (sell) order in the order
book, if such an order exists, and does nothing otherwise.

To formulate this model in more mathematical detail, let $I=(I_-,I_+)\sub\R$
be a nonempty open interval, modeling the possible prices of limit orders,
and let $\ov I:=[I_-,I_+]\sub[-\infty,\infty]$ denote its closure. Let
$\la_\pm:\ov I\to\half$ be functions such that:
\begin{itemize}
\item[(A1)] $\la_-$ is nonincreasing and left-continuous, while $\la_+$ is
  nondecreasing and right-continuous.
\item[(A2)] $\lim_{x\down I_-}\la_-(x)=\la_-(I_-)$ and $\lim_{x\up
  I_+}\la_+(x)=\la_+(I_+)$.
\end{itemize}
We interpret $\la_-(x)$ and $\la_+(x)$ as the \emph{demand} and \emph{supply}
functions, which describe how many items per time unit traders are willing to
buy or sell at the price level $x$. More precisely, let $\mu_\pm$ be finite
measures on $\ov I$ such that
\be\label{mupm}
\mu_-\big([x,I_+]\big)=\la_-(x)
\quand
\mu_+\big([I_-,x]\big)=\la_+(x)
\qquad(x\in\ov I).
\ee
Then the restriction of $\mu_-$ (resp.\ $\mu_+$) to $I$ will be the Poisson
intensity at which traders place buy (resp.\ sell) limit orders at a given
price, while $\mu_-(\{I_+\})$ (resp.\ $\mu_+(\{I_-\})$) will be the Poisson
intensity at which traders place buy (resp.\ sell) market orders. Note that
$\mu_-(\{I_-\})=0=\mu_+(\{I_+\})$ by assumption (A2).

We let $\tau_k$ $(k\geq 1)$ denote the time when the $k$-th trader arrives at
the market, we let $\sig_k\in\{-,+\}$ be a random variable that indicates
whether this trader wants to buy $(-)$ or sell $(+)$, and we let $U_k\in\ov I$
denote the price associated with this trader, where $U_k\in I$ for limit
orders and $U_k=I_\pm$ for market orders. Then
\be\label{Pidef}
\Pi=\{(U_k,\sig_k,\tau_k):k=1,2,\ldots\}
\quad\mbox{with}\quad0<\tau_1<\tau_2<\cdots
\ee
is a Poisson point process on $\ov I\times\{-,+\}\times\half$ with intensity
$\mu\otimes\ell$, where $\ell$ is the Lebesgue measure on $\half$ and $\mu$ is
the finite measure on $\ov I\times\{-,+\}$ given by
$\mu\big(\{\sig\}\times A\big)=\mu_\sig(A)$
for all $\sig\in\{-,+\}$ and measurable $A\sub\ov I$. We let
\be
|\mu_\pm|:=\mu_\pm(\ov I)
\quand
|\mu|:=\mu\big(\ov I\times\{-,+\}\big)=|\mu_-|+|\mu_+|
\ee
denote the total masses of the measures $\mu_\pm$ and $\mu$. To avoid
trivialities, we assume that $|\mu|\neq 0$. Our assumption
that the point process $\Pi$ in (\ref{Pidef}) is Poisson with intensity
$\mu\otimes\ell$ implies that $(\tau_k-\tau_{k-1})_{k\geq 1}$ are
i.i.d.\ exponentially distributed with mean $1/|\mu|$. Moreover, the random
variables $(U_k,\sig_k)_{k\geq 1}$ are i.i.d.\ with law
$\ov\mu:=|\mu|^{-1}\mu$ and independent of $(\tau_k)_{k\geq 1}$.

By definition, a \emph{counting measure} is a measure that can be written as
a finite or countable sum of delta measures. We represent the state of the
order book at a time $t\geq 0$ by a signed counting measure of the form
\be\label{signed}
\Xc=\Xc^+-\Xc^-,
\ee
where $\Xc^-$ and $\Xc^+$ are counting measures on $I$ satisfying
\be\ba{r@{\ }l}\label{signed2}
{\rm(i)}&\mbox{there are no $x,y\in I$ such that $x\leq y$, $\Xc^+(\{x\})>0$,
  $\Xc^-(\{y\})>0$,}\\[5pt]
{\rm(ii)}&\Xc^-\big((x,I_+)\big)<\infty\mbox{ and }
  \Xc^+\big((I_-,x)\big)<\infty\mbox{ for all }x\in I.
\ec
We interpret $\Xc^-(A)$ (resp.\ $\Xc^+(A)$) as the number of buy (resp.\ sell)
limit orders in a measurable set $A\sub I$, and let $\Si_{\rm ord}$ denote the
set of all signed measures that can be written in the form (\ref{signed}) with
$\Xc^\pm$ satisfying (\ref{signed2}). For any $\Xc\in\Si_{\rm ord}$, we let
\bc\label{Mpm}
\dis M_-(\Xc)&:=&\dis\max\big(\{I_-\}\cup\{x\in I:\Xc(\{x\})<0\}\big),\\[5pt]
\dis M_+(\Xc)&:=&\dis\min\big(\{I_+\}\cup\{x\in I:\Xc(\{x\})>0\}\big),
\ec
which can be interpreted as the highest bid and lowest ask price in the order
book.

The state of our Markov process changes only at the times $\tau_1,\tau_2,\ldots$
and we denote the corresponding embedded Markov chain by
\be\label{embed}
X_k:=\Xc_{\tau_k}\qquad(k\geq 0)
\qquad\mbox{with}\quad\tau_0:=0.
\ee
Our previous informal description of the model then translates into the
following definition. Given the initial state $X_0\in\Si_{\rm ord}$, we
inductively define $(X_k)_{k\geq 1}$ as
\be\label{rules}
X_k:=L_{U_k,\sig_k}(X_{k-1})\qquad(k\geq 1),
\ee
where for each $(u,\sig)\in\ov I\times\{-,+\}$, we define a
``Luckock map'' $L_{u,\sig}:\Si_{\rm ord}\to\Si_{\rm ord}$ by
\be\label{Luckmap}
L_{u,\sig}(\Xc):=\left\{\ba{ll}
\dis\Xc-\de_{u\wedge M_+(\Xc)}\quad
&\mbox{if }\sig=-,\ u\wedge M_+(\Xc)\in I,\\[5pt]
\dis\Xc+\de_{u\vee M_-(\Xc)}\quad
&\mbox{if }\sig=+,\ u\vee M_-(\Xc)\in I,\\[5pt]
\dis\Xc\quad&\mbox{otherwise.}
\ea\right.
\ee
For example, if $\sig=+$, then this says that a new sell limit order is added
at the price $u$, unless the current best buy offer $M_-(\Xc)$ is higher than
$u$, in which case this offer is taken, which amounts to adding a delta
measure at $M_-(\Xc)$. The rules for sell market orders are the same, except
that these are not added to the order book if no suitable buy offer exists.

It is easy to see that $(X_k)_{k\geq 0}$ is a Markov chain; in fact, using
terminology from \cite{LPW09}, we have just given a \emph{random mapping
  representation} for it. We call the Markov chain in (\ref{rules}) or, more
or less equivalently, the corresponding continuous-time Markov process
$(\Xc_t)_{t\geq 0}$ the \emph{Stigler-Luckock model} with parameters $\la_\pm$.
In the special case that there are no market orders, this is the model
introduced in \cite{Luc03}. The authors \cite{Sti64,Pla11} considered only the
case that $\mu_-=\mu_+$ is the uniform distribution on a set of 10, resp.\ 100
prices. In \cite{KY16}, infinite piles of limit orders are used as a
construction that is mathematically equivalent to market orders. As we will
see, the introduction of market orders is natural also from a mathematical
point of view and helps us understand the model without market orders.

\begin{figure}[t]
\begin{center}
\begin{tikzpicture}[scale=1]
\begin{scope}[scale=5]
\draw[very thick] (0,0) -- (1,0);
\draw[very thick] (0,0.03) -- (0,-0.03) node[below] {0};
\draw[very thick] (1,0.03) -- (1,-0.03) node[below] {1};

\draw (0.07799298240616463,0) -- ++ (0,-.2);
\draw (0.06843829846385067,0) -- ++ (0,-.2);
\draw (0.06386441502484301,0) -- ++ (0,-.2);
\draw (0.0567212270672393,0) -- ++ (0,-.2);
\draw (0.05200260472267342,0) -- ++ (0,-.2);
\draw (0.04345944334646905,0) -- ++ (0,-.2);
\draw (0.03791787812935048,0) -- ++ (0,-.2);
\draw (0.03584869557543151,0) -- ++ (0,-.2);
\draw (0.0315474254094779,0) -- ++ (0,-.2);
\draw (0.02953613521800729,0) -- ++ (0,-.2);
\draw (0.01823399852607048,0) -- ++ (0,-.2);
\draw (0.01678881351493582,0) -- ++ (0,-.2);
\draw (0.01324015053136334,0) -- ++ (0,-.2);
\draw (0.01063415321503685,0) -- ++ (0,-.2);
\draw (0.008538036324027099,0) -- ++ (0,-.2);
\draw (0.002769627541160657,0) -- ++ (0,-.2);
\draw (0.002134987800456778,0) -- ++ (0,-.2);

\draw (0.4745115843237469,0) -- ++ (0,.2);
\draw (0.5773365771609269,0) -- ++ (0,.2);
\draw (0.5861356306999157,0) -- ++ (0,.2);
\draw (0.5892050856965075,0) -- ++ (0,.2);
\draw (0.6026832237506099,0) -- ++ (0,.2);
\draw (0.6855569878594274,0) -- ++ (0,.2);
\draw (0.6957728844411025,0) -- ++ (0,.2);
\draw (0.7143493089246139,0) -- ++ (0,.2);
\draw (0.7249989077193826,0) -- ++ (0,.2);
\draw (0.7395704716799956,0) -- ++ (0,.2);
\draw (0.741438318706122,0) -- ++ (0,.2);
\draw (0.7658680443598119,0) -- ++ (0,.2);
\draw (0.7766011354920639,0) -- ++ (0,.2);
\draw (0.7798417298190804,0) -- ++ (0,.2);
\draw (0.794302396816905,0) -- ++ (0,.2);
\draw (0.8312501693178945,0) -- ++ (0,.2);
\draw (0.8510234750461139,0) -- ++ (0,.2);
\draw (0.8536212992303549,0) -- ++ (0,.2);
\draw (0.8647098053639197,0) -- ++ (0,.2);
\draw (0.8796656088875631,0) -- ++ (0,.2);
\draw (0.8813223364249116,0) -- ++ (0,.2);
\draw (0.8866967360976634,0) -- ++ (0,.2);
\draw (0.8890849569865692,0) -- ++ (0,.2);
\draw (0.8902579730201873,0) -- ++ (0,.2);
\draw (0.8906254649103109,0) -- ++ (0,.2);
\draw (0.8924004478786421,0) -- ++ (0,.2);
\draw (0.8964286509816569,0) -- ++ (0,.2);
\draw (0.9122922340821561,0) -- ++ (0,.2);
\draw (0.9235749013992267,0) -- ++ (0,.2);
\draw (0.929273673068916,0) -- ++ (0,.2);
\draw (0.930109476770236,0) -- ++ (0,.2);
\draw (0.9356297801704985,0) -- ++ (0,.2);
\draw (0.9361759208262296,0) -- ++ (0,.2);
\draw (0.9421371397941976,0) -- ++ (0,.2);
\draw (0.9692631169139492,0) -- ++ (0,.2);
\draw (0.9736182861925213,0) -- ++ (0,.2);
\draw (0.9843739826928187,0) -- ++ (0,.2);
\end{scope}

\begin{scope}[xshift=7cm,yshift=1cm,xscale=5,yscale=20]
\draw[very thick] (0,0) -- (1,0);
\foreach \x in {0.2,0.4,0.6,0.8,1}
 \draw[very thick] (\x,0) -- (\x,-.01) node[below] {\x};
\draw[very thick] (0,-0.11) -- (0,0.03);
\draw[very thick] (0,-0.1) --++ (-.03,0) node[left] {-1000};
\draw[very thick] (0,-0.075) --++ (-.03,0) node[left] {-750};
\draw[very thick] (0,-0.05) --++ (-.03,0) node[left] {-500};
\draw[very thick] (0,-0.025) --++ (-.03,0) node[left] {-250};
\draw[very thick] (0,0) --++ (-.03,0) node[left] {0};
\draw[very thick] (0,0.025) --++ (-.03,0) node[left] {250};

\draw plot file {distrib10000.dat};
\end{scope}

\end{tikzpicture}
\caption{Simulation of the ``uniform'' Stigler-Luckock model with $I=(0,1)$,
  $\la_-(x)=1-x$, and $\la_+(x)=x$. On the left: the state $X_{250}$ of the
  order book after the arrival of 250 traders, starting from the empty initial
  state. On the right: the distribution function $x\mapsto X_k([0,x])$ of the
  random signed measure $X_k$ after the arrival of $k=10,000$ traders.}
\label{fig:distrib}
\end{center}
\end{figure}
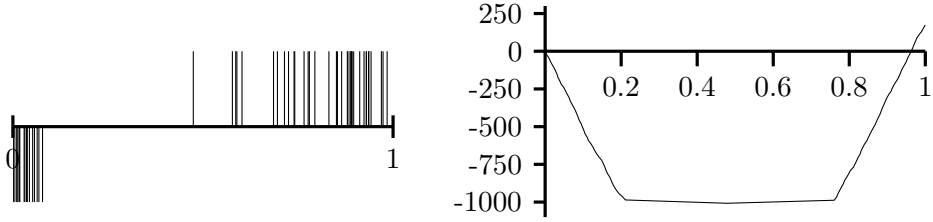

We let
\be
V_{\rm max}:=\la_-(I_-)\wedge\la_+(I_+)
\quand
V_{\rm W}:=\sup_{x\in\ov I}\la_-(x)\wedge\la_+(x)
\ee
denote the maximal possible volume of trade and the maximal possible volume of
trade at a fixed price level, respectively. We will sometimes need the
following stronger conditions on our demand and supply functions.
\begin{itemize}
\item[(A3)] $\la_\pm$ are continuous on $\ov I$, $\la_-$ is nonincreasing, and
  $\la_+$ is nondecreasing.
\item[(A4)] The function $\la_+-\la_-$ is strictly increasing on $\ov I$.
\item[(A5)] The functions $\la_-$ and $\la_+$ are $>0$ on $I$.
\item[(A6)] The rates $\la_+(I_-)$ and $\la_-(I_+)$ of market orders are both
  $>0$.
\item[(A7)] $V_{\rm W}<V_{\rm max}$.
\end{itemize}
In particular, (A3) implies (A1) and (A2). As shown in
Appendix~\ref{A:staform}, (A3) and (A4) are not really a restriction, since
every Stigler-Luckock model satisfying (A1) and (A2) can be obtained as a
function of a Stigler-Luckock model satisfying (A3) and, under mild extra
assumptions, also (A4). Condition~(A5) also comes basically without loss of
generality, since sell orders on the right of the first point $x$ where
$\la_-(x)=0$ are trivially never matched, and similarly for buy orders at the
other end of the interval. The conditions (A6) and (A7) are restrictive, of
course.

Conditions (A3), (A4), and (A7) imply that there exists a unique point $x_{\rm
  W}\in I$ such that
\be\label{xWalras}
\la_-(x_{\rm W})=\la_+(x_{\rm W}).
\ee
Classical economic theory going back to Walras \cite{Wal74} says that in an
infinitely liquid market, the equilibrium price is $x_{\rm W}$, which is why
we call $x_{\rm W}$ the \emph{Walrasian price}. Note that $V_{\rm
  W}:=\la_-(x_{\rm W})=\la_+(x_{\rm W})$, which is why we call this the
\emph{Walrasian volume of trade}.


\subsection{Luckcock's differential equation}\label{S:Luckdif}

The following theorem is essentially proved in \cite{Luc03}, but for
completeness we will provide a proof in the present setting. Below, if $f:\ov
I\to\R$ is a continuous function of bounded variation, then we let $\di f$
denote the signed measure on $\ov I$ such that $\di
f\big((x,y]\big):=f(y)-f(x)$. If $g:\ov I\to\R$ is a bounded measurable
function, then $g\,\di f$ denotes the measure $\di f$ weighted with $g$, i.e.,
$g\,\di f\big((x,y]\big):=\int_x^yg\,\di f$. We call any pair $(f_-,f_+)$ of
  continuous functions of bounded variation such that (\ref{Luckock}) below
  holds a \emph{solution to Luckock's equation}.

\bt[Luckock's differential equation]
Consider\label{T:Luckock} a Stigler-Luck\-ock model with supply and demand
functions $\la_\pm:\ov I\to\half$ satisfying (A3) and (A5). Assume that the
model has an invariant law on $\Si_{\rm ord}$ and let $(X_k)_{k\geq
  0}$ denote the corresponding stationary process. Then the functions
$f_\pm:\ov I\to\R$ defined by
\be\label{fpm}
f_-(x):=\P\big[M_-(X_k)\leq x\big]
\quand
f_+(x):=\P\big[M_+(X_k)\geq x\big]\qquad(x\in\ov I)
\ee
are continuous and solve the equations
\be\ba{rr@{\,}c@{\,}l}\label{Luckock}
{\rm(i)}&\dis f_-\di\la_++\la_-\di f_+&=&\dis0,\\[5pt]
{\rm(ii)}&\dis f_+\di\la_-+\la_+\di f_-&=&\dis0,\\[5pt]
{\rm(iii)}&\multicolumn{3}{c}{\dis f_-(I_+)=1=f_+(I_-).}
\ec
\et

We remark that although Theorem~\ref{T:Luckock} shows that the equilibrium
distributions of the best buy and sell order in the order book can more or
less be solved explicitly (depending on our ability to solve (\ref{Luckock})),
this does not automatically mean that Stigler-Luckock models as a whole are
``solvable''. For example, we do not know how to explicitly calculate the
joint distribution of $M_-(\Xc)$ and $M_+(\Xc)$ (as opposed to its marginals).
Also, it seems to be quite hard to get information about the equilibrium
distribution of seemingly simple functions of the process like the number of
sell (or buy) limit orders in a certain interval.

Theorem~\ref{T:Luckock} motivates the study of solutions to Luckock's equation
(\ref{Luckock}). 

\bp[Solutions to Luckock's equation]
Assume\label{P:Luck} (A3) and (A6). Then Luckock's equation has a unique
solution $(f_-,f_+)$. One has
\be\ba{rrcl}\label{valid}
{\rm(i)}&\dis f_-(I_-)\geq0&\desd&\dis\La_-:=\frac{1}{\la_-(I_-)\la_-(I_+)}
-\int_{I_-}^{I_+}\frac{1}{\la_+}\di\Big(\frac{1}{\la_-}\Big)\geq0,\\[5pt]
{\rm(ii)}&\dis f_+(I_+)\geq0&\desd&\dis\La_+:=\frac{1}{\la_+(I_-)\la_+(I_+)}
+\int_{I_-}^{I_+}\frac{1}{\la_-}\di\Big(\frac{1}{\la_+}\Big)\geq0.
\ec
Both formulas also hold with the inequality signs reversed.
The functions $(f_-,f_+)$ also satisfy
\be\label{volume}
\la_+(I_+)-f_-(I_-)\la_+(I_-)=\la_-(I_-)-f_+(I_+)\la_-(I_+).
\ee
\ep

If the solution $(f_-,f_+)$ to Luckock's equation satisfies $f_-(I_-)\wedge
f_+(I_+)\geq 0$, then we call such a solution \emph{valid}. See
Figure~\ref{fig:numer} for a plot of $(f_-,f_+)$ for one particular model -in
this particular example, $(f_-,f_+)$ is not valid. By Theorem~\ref{T:Luckock},
a necessary condition for a Stigler-Luckock model to have an invariant law is
that the solution to Luckock's equation is valid. We conjecture that this
condition is also sufficient, but stop short of proving this. (Note however
Theorem~\ref{T:posrec} below, which goes some way in this direction.)

If a Stigler-Luckock model has an invariant law, then the quantity in
(\ref{volume}) can be interpreted as the \emph{volume of trade}, i.e., the
expected number of orders (of either type) that are matched per unit of
time. Indeed, since the process has an invariant law, sell limit orders, which
arrive at rate $\la_+(I_+)-\la_+(I_-)$, are all eventually matched, while
$1-f_-(I_-)=\P[M_-(\Xc_t)>I_-]$ is the fraction of sell market orders that are
matched, so $\big(\la_+(I_+)-\la_+(I_-)\big)+\big(1-f_-(I_-)\big)\la_+(I_-)$
is the total rate at which sell orders are matched, which equals the left-hand
side of (\ref{volume}). The right-hand side of (\ref{volume}) has a similar
interpretation in terms of buy orders.


\subsection{Positive recurrence}

Let $(X_k)_{k\geq 0}$ be a Stigler-Luckock model with discrete time (i.e., the
embedded Markov chain from (\ref{embed})), started in the empty initial state
$X_0=0$, and let $\tau$ denote its first return time to $0$, i.e.,
$\tau:=\inf\{k>0:X_k=0\}$. We say that a Stigler-Luckock model is
\emph{positive recurrent} if $\E[\tau]<\infty$, \emph{transient} if
$\P[\tau=\infty]>0$, and \emph{null recurrent} in the remaining case.  The
main result of the present paper is the following result, that gives a more or
less complete characterization of positive recurrent Stigler-Luckock
models. Below and in what follows, we let $\Si^{\rm fin}_{\rm ord}$ denote the
set of all finite configurations $\Xc\in\Si_{\rm ord}$, i.e., those for which
$\Xc^-$ and $\Xc^+$ are finite measures.

\bt[Positive recurrence]
Assume\label{T:posrec} (A3) and (A6). Then a Stig\-ler-Luckock model is positive recurrent if
and only if the unique solution $(f_-,f_+)$ to Luckock's equation satisfies
$f_-(I_-)\wedge f_+(I_+)>0$. If a Stigler-Luckock model is positive recurrent,
then it has an invariant law $\nu$ that is concentated on $\Si^{\rm fin}_{\rm
  ord}$. Moreover, the process started in any initial law that is concentated
on $\Si^{\rm fin}_{\rm ord}$ satisfies
\be\label{tonu}
\big\|\P[X_k\in\,\cdot\,]-\nu\big\|\asto{k}0,
\ee
where $\|\,\cdot\,\|$ denotes the total variation norm.
\et

Under rather restrictive additional assumptions, a version of
Theorem~\ref{T:posrec} has also been proved in \cite[Thm~2.1~(3)]{KY16}.
We will prove Theorem~\ref{T:posrec} by constructing a Lyapunov
function. Instead, the authors of \cite{KY16}, being unable to find a Lyapunov
function, used fluid limit methods.

\subsection{Restricted models}\label{S:restrict}

Assume that the demand and supply functions $\la_-$ and $\la_+$ satisfy (A3)
and (A5). Then, for each interval $J=(J_-,J_+)$ such that $\ov J=[J_-,J_+]\sub
I$, the restrictions of $\la_-$ and $\la_+$ to $\ov J$ satisfy (A6). We call
the corresponding Stigler-Luckock model the \emph{restricted model on $J$}.
By Proposition~\ref{P:Luck}, Luckock's equation has a unique solution for this
restricted model, and by Theorem~\ref{T:posrec} we can read off from this
solution whether the restricted model is positive recurrent. In the present
section, for fixed $I$ and $\la_\pm$, we investigate the set of all
subintervals $\ov J\sub I$ for which the restricted model is positive
recurrent.

We note that if $(X_k)_{k\geq 0}$ is a Stigler-Luckock model on $I$ and
$X_k\big|_J$ denotes the restriction of the random signed measure $X_k$ to a
subinterval $J\sub I$, then it is in general not true that
$(X_k\big|_J)_{k\geq 0}$ is a Markov chain. In particular, this is not the
same as the restricted model on $J$. Nevertheless, we will see that under
suitable conditions, there exists a special subinterval $J\sub I$ (below, this
is called the \emph{competitive window}) such that in the long run, we expect
the model on $I$ to behave basically like the model restricted to $J$, with
all buy limit orders on the left of $J$ and all sell limit orders on the right
of $J$ never being matched and as a result staying in the order book forever.

Assume (A3) and (A5) and for $I_-<J_-<J_+<I_+$, let $\La_-(J_-,J_+)$ and
$\La_+(J_-,J_+)$ denote the expressions in (\ref{valid}), calculated for the
process restricted to the subinterval $\ov J$. For fixed $J_-\in I$
resp.\ $J_+\in I$, we define
\bc\label{phipm}
\dis\phi_-(J_+)
&:=&\dis\sup\big\{J_-\in(I_-,J_+):\La_-(J_-,J_+)\leq 0\big\},\\[5pt]
\dis\phi_+(J_-)&:=&\dis\inf\big\{J_+\in(J_-,I_+):\La_+(J_-,J_+)\leq 0\big\},
\ec
with the conventions $\sup\emptyset:=I_-$ and $\inf\emptyset:=I_+$.
Let
\be\ba{r@{\,}l}\label{Rdef}
\dis R:=\big\{(J_-,J_+)\in I\times I:
&\dis J_-<J_+\mbox{ and the restricted}\\
&\dis\mbox{model on $J$ is positive recurrent}\big\}.
\ec
The following lemma says that the set $R$ is bounded by the graphs of the
functions $\phi_\pm$, as well as (trivially) the line $J_-=J_+$.

\bl[Positive recurrence of restricted models]\label{L:JJ}
Assume (A3) and (A5). Then $\phi_-(J_+)<J_+$ and $J_-<\phi_+(J_-)$ for all
$J_-,J_+\in I$. Moreover, a point $(J_-,J_+)\in I\times I$ belongs to the set
$R$ from (\ref{Rdef}) if and only $\phi_-(J_+)<J_-$, $J_+<\phi_+(J_-)$, and
$J_-<J_+$.
\el

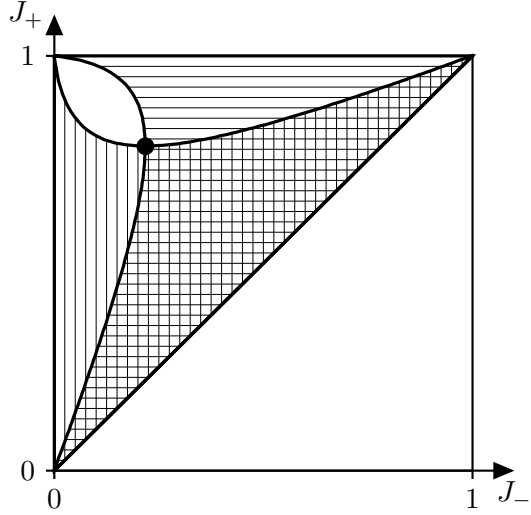
\begin{figure}[t]
\begin{center}
\begin{tikzpicture}[scale=5.5,>=triangle 45]
\draw[thick,->] (0,0)--(1.1,0) node[below] {$J_-$}; 
\draw[thick,->] (0,0)--(0,1.1) node[left] {$J_+$};
\draw[thick] (0,0)--++(0,-0.02) node[below] {$0$};
\draw[thick] (1,0)--++(0,-0.02) node[below] {$1$};
\draw[thick] (0,0)--++(-0.02,0) node[left] {$0$};
\draw[thick] (0,1)--++(-0.02,0) node[left] {$1$};
\draw[thick] (0,0)--(1,0)--(1,1)--(0,1)--cycle;
\draw[very thick] (0,0)--(1,1);
\draw[very thick] plot file {phip.dat};
\draw[very thick] plot file {phim.dat};
\filldraw (0.2178,0.7822) circle (0.02cm);
\draw[clip] (0,0)--(1,1)--(0,1)--cycle;
\begin{scope}
\draw[clip] plot file {phip.dat};
\foreach \x in {0,0.025,...,1}
 \draw[very thin] (\x,0) -- (\x,1); 
\end{scope}
\begin{scope}
\draw[clip] plot file {phim.dat};
\foreach \x in {0,0.025,...,1}
 \draw[very thin] (0,\x) -- (1,\x); 
\end{scope}
\end{tikzpicture}
\caption{Restrictions of the uniform Stigler-Luckock model with $I=(0,1)$,
  $\la_-(x)=1-x$, and $\la_+(x)=x$ to a subinterval $(J_-,J_+)$. The solution
  $(f_-,f_+)$ to Luckock's equation for the restricted model satisfies
  $f_+(J_+)>0$ in the vertically striped area and $f_-(J_-)>0$ in the
  horizontally striped area. The intersection of these areas corresponds to
  the set $R$ of subintervals on which the restricted model is positive
  recurrent. The intersection of the curves $J_-=\phi_-(J_+)$ and
  $J_+=\phi_+(J_-)$, indicated with a dot, corresponds to the competitive
  window.}
\label{fig:R}
\end{center}
\end{figure}

In Figure~\ref{fig:R} we have pictured the set $R$ and the graphs of the
functions $\phi_\pm$ for the ``uniform'' model with $I=[0,1]$, $\la_-(x)=1-x$,
and $\la_+(x)=x$. For this model, one can check that the solution of Luckock's
equation for the restricted model on $J$ satisfies $f_-(J_-)=0$ if and only if
$J_-=\phi_-(J_+)$, and likewise one has $f_+(J_+)=0$ if and only if
$(J_-,J_+)$ lies on the graph $\{J_+=\phi_+(J_-)\}$. The graphs of the
functions $\phi_\pm$ intersect in a unique point, which in the light of
(\ref{volume}) must satisfy $\la_-(J_-)=\la_+(J_+)$.

\subsection{The competitive window}\label{S:compwind}

In this subsection, we focus on subintervals $J=(J_-,J_+)$ that are
``symmetric'' in the sense that $\la_-(J_-)=\la_+(J_+)$. Let
$\la^{-1}_-:[\la_-(I_+),\la_-(I_-)]\to\ov I$ and
$\la^{-1}_+:[\la_+(I_-),\la_+(I_+)]\to\ov I$ denote the left-continuous
inverses of the functions $\la_-$ and $\la_+$, respectively, i.e.,
\bc\label{jpm}
\dis\la^{-1}_-(V)&:=&\dis\sup\{x\in\ov I:\la_-(x)\geq V\},\\[5pt]
\dis\la^{-1}_+(V)&:=&\dis\inf\{x\in\ov I:\la_+(x)\geq V\}.
\ec
Assuming (A5) and (A7), we define a function $\Phi:[V_{\rm W},V_{\rm
    max}]\to[0,\infty]$ by
\be\label{Phidef}
\Phi(V):=-\int_{V_{\rm W}}^V\Big\{\frac{1}{\la_+\big(\la^{-1}_-(W)\big)}
+\frac{1}{\la_-\big(\la^{-1}_+(W)\big)}\Big\}\di\Big(\frac{1}{W}\Big).
\ee
Note that $\Phi$ is increasing since $1/W$ is a decreasing function.

\bp[Symmetric subintervals]\hspace*{1cm}
Assume\label{P:Phi} (A3), (A5), and (A7), let $V\in(V_{\rm W},V_{\rm max}]$,
assume that $J=(J_-,J_+):=\big(\la^{-1}_-(V),\la^{-1}_+(V)\big)$ satisfies
$\ov J\sub I$, and let $(f_-,f_+)$ be the unique solution to Luckock's
equation on $J$. Then:
\begin{enumerate}
\item If $\Phi(V)<V_{\rm W}^{-2}$, then $f_-(J_-)>0$ and $f_+(J_+)>0$.
\item If $\Phi(V)=V_{\rm W}^{-2}$, then $f_-(J_-)=0$ and $f_+(J_+)=0$.
\item If $\Phi(V)>V_{\rm W}^{-2}$, then $f_-(J_-)<0$ and $f_+(J_+)<0$.
\end{enumerate}
\ep

Theorem~\ref{T:posrec} and Proposition~\ref{P:Phi} motivate us to define
\emph{Luckock's volume of trade} as
\be\label{VL}
V_{\rm L}:=\sup\big\{V\in[V_{\rm W},V_{\rm max}]:\Phi(V)\leq V_{\rm W}^{-2}\big\},
\ee
and the \emph{competitive window} as $J^{\rm c}=(J^{\rm c}_-,J^{\rm
  c}_+):=\big(\la^{-1}_-(V_{\rm L}),\la^{-1}_+(V_{\rm L})\big)$. Further
motivation for this definition comes from the following observation.

\bl[Competitive window]
Assume\label{L:comp} (A3), (A5), and (A7), and that $V_{\rm L}<V_{\rm max}$.
Then there exists a unique subinterval $J$ such that $\ov J\sub I$ and:
\begin{enumerate}
\item The unique solution $(f_-,f_+)$ to Luckock's equation on $J$ satisfies\\
  $f_-(J_-)=0=f_+(J_+)$.
\item $\la_-<\la_-(J_-)$ on $(J_-,J_+]$ and $\la_+<\la_+(J_+)$ on $[J_-,J_+)$.
\end{enumerate}
This subinterval is given by $J=J^{\rm c}$.
\el

We note that $V_{\rm W}<V_{\rm L}$ always, but it is possible the $V_{\rm
  L}=V_{\rm max}$. For example, this happens for the model with $\ov I=[0,1]$,
$\la_-(x)=(1-x)^\al$, and $\la_+(x)=x^\al$ if $0<\al\leq1/2$. Another example
are models of the form $\ov I=[0,1+\la]$, $\la_-(x)=(1+\la-x)\wedge 1$, and
$\la_+(x)=x\wedge 1$, where $\la\in[1-1/z,1)$ with $z$ as in
Lemma~\ref{L:unif} below; compare the discussion in \cite[Section~5.1]{KY16}.

Since $V_{\rm L}$ us usually much larger than $V_{\rm W}$, the Stigler-Luckock
model is highly non-liquid. As such, it is not a realistic model of a real
market, though it may be a useful first step towards building more realistic
models. The special case where buy and sell limit orders are uniformly
distributed on the unit interval is of some special interest. Numerically, the
constant $V_{\rm L}$ from Lemma~\ref{L:unif} is given by
$V_{\rm L}\approx0.78218829428020$.

\bl[Uniform model]\label{L:unif}
The Stigler-Luckock model with $\ov I=[0,1]$, $\la_-(x)=1-x$, and $\la_+(x)=x$
has a competitive window $(J^\cri_-,J^\cri_+)$ which is given by
$1-J^\cri_-=J^\cri_+=V_{\rm L}$, where $V_{\rm L}=1/z$ with $z$ the unique
solution of the equation $e^{-z}-z+1=0$.
\el

\subsection{Discussion and open problems}\label{S:open}

Kelly and Yudovina \cite[Thm~2.1 (1) and (2)]{KY16} have shown\footnote{The
  statement of Thm~2.1~(1) in the preprint \cite{KY16} is not entirely
  correct (and will be corrected in the final publication), but their proof of
  (\ref{Msupinf}) is correct.} that there exists deterministic constants
$J^\ast_\pm$ such that
\be\label{Msupinf}
\liminf_{k\to\infty}M_-(X_k)=J^\ast_-
\quand\limsup_{k\to\infty}M_+(X_k)=J^\ast_+\quad{\rm a.s.}.
\ee
For the uniform model of Lemma~\ref{L:unif}, they moreover show
\cite[Corollary~2.3]{KY16} that $J^\ast_\pm=J^\cri_\pm$, where
$(J^\cri_-,J^\cri_+)$ is the competitive window as defined below
(\ref{VL}). Their methods apply more generally, although they need certain
technical assumptions, in particular, that the measures $\mu_\pm$ from
(\ref{mupm}) have densities w.r.t.\ the Lebesgue measure that are bounded away
from zero and infinity.

We conjecture that $\liminf_{k\to\infty}\la_\pm(M_\pm(X_k))=V_{\rm L}$
a.s.\ holds generally under the assumptions (A3), (A5) and (A7). Assuming
moreover that $V_{\rm L}<V_{\rm max}$, we conjecture that if
$X_k\big|_{J^\cri}$ denotes the restriction of $X_k$ to $J^\cri$, then,
starting from any finite initial state, the law of $X_k\big|_{J^\cri}$
converges as $k\to\infty$ to a limit law that is concentrated on $\Si_{\rm
  ord}$, and that this limit law is an invariant law for the restricted model
on $J^\cri$. Indeed, (\ref{Msupinf}) says that in the long run, the price of
the best buy offer never drops below $J^\cri_-$ and the price of the best sell
offer never climbs above $J^\cri_+$, which allows us to treat limit sell
orders at prices below $J^\cri_-$ and limit buy orders above $J^\cri_+$ as
market orders. Simulations (see Figure~\ref{fig:distrib}) seem to support
these conjectures.

Proving these conjectures remains an open problem. Theorems~\ref{T:Luckock}
and \ref{T:posrec} allow us to conclude, however, that if
$\la_\pm(J_\pm)>V_{\rm L}$, then the restricted model on $J$ does not have an
invariant law while if $\la_\pm(J_\pm)<V_{\rm L}$, then the restricted model
on $J$ is positive recurrent. Further motivation for the conjectures comes
from the study of similar models. In \cite{Swa15}, a ``one-sided canyon
model'' is studied that is in many ways similar to the Stigler-Luckock model
except that there is only one type of points as opposed to the two types (buy
and sell orders) of a Stigler-Luckock model. This ``one-sided'' model also has
a competitive window that can be calculated explicitly and in fact the
analogues of the conjectures above have all been proved for this model, mainly
due to the hugely simplifying fact that for this model, restricting the
process to a smaller interval does again yield a Markov chain.

In this context, we also mention a model for email communication due to
Gabrielli and Caldarelli \cite{CG09}. This model is even simpler than the
previous one since not only is the restriction of the process to a subinterval
Markovian, but even just counting the number of points in a subinterval
already yields a Markov chain. For this model, it has been possible to solve
subtle questions about the behavior of the stationary process near the
boundary of the competitive window \cite{FS16}.

The models mentioned so far belong to a wider class of models that also
includes the Bak Sneppen model \cite{BS93} and its modified version from
\cite{MS12}, as well as the branching Brownian motions with strong selection
treated in \cite{Mai16}. All these models implement some version of the rule
``kill the lowest particle'' and seem to exhibit self-organized criticality,
although this has been rigorously proved only for some of the models.

As mentioned before, the Stigler-Luckock model describes an extremely
non-liquid market, and (mainly) for that reason is not a realistic model for a
real market, although it may perhaps be used as a first step towards more
realistic models. In recent years, there has been considerable activity in the
search for simple, yet realistic models for an order book. An attempt
to make the Stigler-Luckock model more realistic by introducing market makers
is made in \cite{PS16}. Often, authors assume that new orders arrive relative
to the current best bid or ask price
\cite{Mas00,CST10,Kru12,LRS13,LRS16,SRR16}. For mathematical simplicity, it is
sometimes assumed that the order book is ``one-sided'' in the sense that all
buy orders are market orders and all sell orders are limit orders
\cite{LRS13,Tok15,LRS16}. A very general model is formulated in \cite{Smi12}.
We refer the reader to the overview artice \cite{CTPA11} or Chapter~4 of the
book \cite{Sla13}, for a more complete view on this topic.



\subsection{Methods}

The results in Sections~\ref{S:Luckdif}, \ref{S:restrict}, and
  \ref{S:compwind} are mainly a reworking of similar results already proved
by Luckock in \cite{Luc03}, although Proposition~\ref{P:Phi} is a significant
improvement over \cite[Prop.~4]{Luc03}. Nevertheless, Luckock already derived
the differential equation (\ref{Luckock}) and showed how it could be used to
calculate the competitive window for a given model. Throughout his paper,
however, he takes stationarity as a model assumption, where in fact,
``stationarity'' for him means the existence of two prices $J_-<J_+$ such that
buy orders on the left of $J_-$ and sell orders on the right of $J_+$ are
never matched while the process inside $(J_-,J_-)$ is stationary in law.

From a mathematical point of view, the existence of such a stationary process
requires proof. Moreover, one would like to prove that the process started in
an arbitrary finite initial state converges, in a suitable sense, to such a
stationary state. For Luckock's original model, these problems remain open,
but for positive recurrent processes with market orders, these questions are
resolved by Theorem~\ref{T:posrec}, which is the most important contribution
of the present paper.

The proof of Theorem~\ref{T:posrec} is based on a Lyapunov function. We equip
the space $\Si^{\rm fin}_{\rm ord}$ with a topology such that $\Xc(n)\to\Xc$
if and only if $\Xc^\pm(n)\Rightarrow\Xc^\pm$, where $\Rightarrow$ denotes
weak convergence. We also equip $\Si^{\rm fin}_{\rm ord}$ with the Borel
\si-algebra associated with this topology. We call a function $F:\Si^{\rm
  fin}_{\rm ord}\to\R$ \emph{Lipschitz} if there exists a constant $L$ such
that $|F(\Xc+\de_x)-F(\Xc)|\leq L$ for all $x\in I$. For any measurable
Lipschitz function $F:\Si^{\rm fin}_{\rm ord}\to\R$, write
\be\label{Gdef}
GF(\Xc):=\int\big\{F\big(L_{u,\sig}(\Xc)\big)-F\big(\Xc\big)\big\}
\,\mu\big(\di(u,\sig)\big),
\ee
where $L_{u,\sig}$ is the Luckock map defined in (\ref{Luckmap}) and $\mu$ is
the measure defined below (\ref{Pidef}). Then $G$ is the generator of the
continuous-time Markov process $(\Xc_t)_{t\geq 0}$.

It turns out that there is a useful and explicit formula for $GF$ when $F$ is
a ``linear'' function of the form
\be\label{Flin}
F(\Xc):=\int_Iw_-(x)\Xc^-(\di x)+\int_Iw_+(x)\Xc^+(\di x)
\qquad\big(\Xc\in\Si^{\rm fin}_{\rm ord}\big),
\ee
where $w_\pm:\ov I\to\R$ are bounded ``weight'' functions such that $w_-$ is
left-continuous and $w_+$ is right-continuous. The values of $w_-$ and $w_+$
in the boundary points $I_-$ and $I_+$ are irrelevant for (\ref{Flin}), but
for notational convenience, we define $w_-(I_+)$ and $w_+(I_-)$ by left,
resp.\ right continuity, and use the convention that
\be\label{wbc}
w_-(I_-):=0\quand w_+(I_+):=0.
\ee
With this convention, the following lemma describes the action of the
generator on linear functions of the form (\ref{Flin}).

\bl[Generator on linear functionals]
Assume\label{L:GF} (A3). Then, for functions of the form (\ref{Flin}), one has
\be\label{GF}
GF(\Xc)=q_-\big(M_-(\Xc)\big)+q_+\big(M_+(\Xc)\big)
\qquad\big(\Xc\in\Si^{\rm fin}_{\rm ord}\big),
\ee
where $q_\pm:\ov I\to\R$ are given by
\be\ba{r@{\,}c@{\,}l}\label{qpm}
\dis q_-(x)&:=&\dis\int_x^{I_+}w_+\di\la_+-w_-(x)\la_+(x),\\[5pt]
\dis q_+(x)&:=&\dis-\int_{I_-}^xw_-\di\la_--w_+(x)\la_-(x).
\ec
\el
\bpro
We observe that $\int_{M_-(\Xc)}^{I_+}w_+\di\la_+$ is the rate at which
$F(\Xc)$ increases due to sell limit orders being added to the order book
while 
the product of $w_-\big(M_-(\Xc)\big)$ and $\la_+\big(M_-(\Xc)\big)$
is the rate at which
$F(\Xc)$ decreases due to buy limit orders being removed from the order
book. In view of our convention (\ref{wbc}), the latter term is zero when
the order book contains no buy limit orders. The two terms in
$q_+(M_+(\Xc))$ in have similar interpretations.
\epro

Formula (\ref{qpm}) tells us how to calculate the functions $q_\pm$ from
(\ref{GF}) from the weight functions $w_\pm$. It turns out that under the
assumptions (A3) and (A6), one can uniquely solve the following inverse
problem: if $q_\pm$ are given up to an additive constant, then find $w_\pm$
such that (\ref{qpm}) holds. This is shown in Theorem~\ref{T:inverse} below
and more specifically for indicator functions of the form $q_-=1_{[I_-,z]}$
and $q_+=1_{[z,I_+]}$ in the following theorem, that moreover specifies the
additive constant.

\bt[Special weight functions]\label{T:spec}
Assume (A3) and (A6). Then, for each $z\in\ov I$, there exist a unique pair of
bounded weight functions
\[
(w^{z,-}_-,w^{z,-}_+)=(w_-,w_+)
\]
such that $w_-$ is
left-continuous and $w_+$ is right-continuous, and the linear
functional $F^{z,-}=F$ from (\ref{Flin}) satisfies
\be\label{GFfm}
GF(\Xc)=1_{\txt\{M_-(\Xc)\leq z\}}-f_-(z)\qquad\big(\Xc\in\Si^{\rm fin}_{\rm ord}\big),
\ee
where $(f_-,f_+)$ is the unique solution to Luckock's equation
(\ref{Luckock}). Likewise, there exist a unique pair of weight functions
$(w^{z,+}_-,w^{z,+}_+)=(w_-,w_+)$ such that the linear functional $F^{z,+}=F$
from (\ref{Flin}) satisfies
\be\label{GFfp}
GF(\Xc)=1_{\txt\{M_+(\Xc)\geq z\}}-f_+(z)\qquad\big(\Xc\in\Si^{\rm fin}_{\rm ord}\big),
\ee
\et

\begin{figure}[t]
\begin{tikzpicture}[>=triangle 45,scale=0.8]
\begin{scope}[xscale=6,yscale=4]
\draw[->] (0,0) -- (1.1,0);
\draw[->] (0,0) -- (0,1);
\draw[thick] (-0.025,0) node[left] {0};
\draw[thick] (0.025,0.3) -- (-0.025,0.3) node[left] {0.3};
\draw[thick] (0.025,0.6) -- (-0.025,0.6) node[left] {0.6};
\draw[thick] (0.025,0.9) -- (-0.025,0.9) node[left] {0.9};
\draw[thick] (0,-0.025) node[below] {0};
\draw[thick] (0.3,0.025) -- (0.3,-0.025) node[below] {0.3};
\draw[thick] (0.9,0.025) -- (0.9,-0.025) node[below] {0.9};
\draw[blue,dashed,very thick] plot file {lam.dat};
\draw[red,very thick] plot file {lap.dat};
\draw (0.35,0.7) node[above,blue] {$\la_-$};
\draw (0.8,0.75) node[below,red] {$\la_+$};
\end{scope}
\begin{scope}[xscale=6,yscale=1.7,xshift=1.3cm,yshift=1cm]
\draw[->] (-0.025,0) node[left] {0} -- (1.1,0);
\draw[thick] (0.025,1) -- (-0.025,1) node[left] {1};
\draw[thick] (0.025,0.5) -- (-0.025,0.5) node[left] {0.5};
\draw[thick] (0.025,-0.5) -- (-0.025,-0.5) node[left] {-0.5};
\draw[thick] (0.025,-1) -- (-0.025,-1) node[left] {-1};
\draw[->] (0,-1.1) -- (0,1.4);
\draw[thick] (0.3,0.05) -- (0.3,-0.05) node[below] {0.3};
\draw[thick] (0.9,0.05) -- (0.9,-0.05) node[below] {0.9};
\draw[blue,dashed,very thick] plot file {fm.dat};
\draw[red,very thick] plot file {fp.dat};
\draw (0.2,1) node[red] {$f_+$};
\draw (0.8,0.8) node[blue] {$f_-$};
\end{scope}
\begin{scope}[xscale=6,yscale=0.35,yshift=-11cm]
\draw[->] (-0.025,0) node[left] {0} -- (1.1,0);
\draw[thick] (0.025,6) -- (-0.025,6) node[left] {6};
\draw[thick] (0.025,3) -- (-0.025,3) node[left] {3};
\draw[thick] (0.025,-3) -- (-0.025,-3) node[left] {-3};
\draw[->] (0,-4.2) -- (0,7.5);
\draw[thick] (0.3,0.25) -- (0.3,-0.25) node[below] {0.3};
\draw[thick] (0.9,0.25) -- (0.9,-0.25) node[below] {0.9};
\draw[blue,dashed,very thick] plot file {wMm.dat};
\draw[red,very thick] plot file {wMp.dat};
\draw (0.4,3) node[blue,above] {$w^{(-)}_-$};
\draw (0.6,-2) node[red,below] {$w^{(-)}_+$};
\end{scope}
\begin{scope}[xscale=6,yscale=0.35,xshift=1.3cm,yshift=-11cm]
\draw[->] (-0.025,0) node[left] {0} -- (1.1,0);
\draw[thick] (0.025,6) -- (-0.025,6) node[left] {6};
\draw[thick] (0.025,3) -- (-0.025,3) node[left] {3};
\draw[thick] (0.025,-3) -- (-0.025,-3) node[left] {-3};
\draw[->] (0,-4.2) -- (0,7.5);
\draw[thick] (0.3,0.25) -- (0.3,-0.25) node[below] {0.3};
\draw[thick] (0.75,0.25) -- (0.75,-0.25) node[below] {$z$};
\draw[thick] (0.9,0.25) -- (0.9,-0.25) node[below] {0.9};
\draw[blue,dashed,very thick] plot file {wm.dat};
\draw[red,very thick] plot file {wp.dat};
\draw (0.5,-1.8) node[blue,below] {$w^{z,+}_-$};
\draw (0.6,5) node[red,above] {$w^{z,+}_+$};
\end{scope}
\end{tikzpicture}
\caption{The solution $(f_-,f_+)$ to Luckock's
  equation, as well as two examples of weight functions $(w_-,w_+)$ as in
  Theorem~\ref{T:spec}. In this example, $I=(0.3,0.9)$,
  $\la_-(x)=1-x$, and $\la_+(x)=x$. The lower left picture shows the weight
  functions $w^{(-)}_\pm=w^{I_-,-}_\pm$ while the lower right picture shows
  the weight functions $w^{z,+}_\pm$ for $z=0.75$.}
\label{fig:numer}
\end{figure}
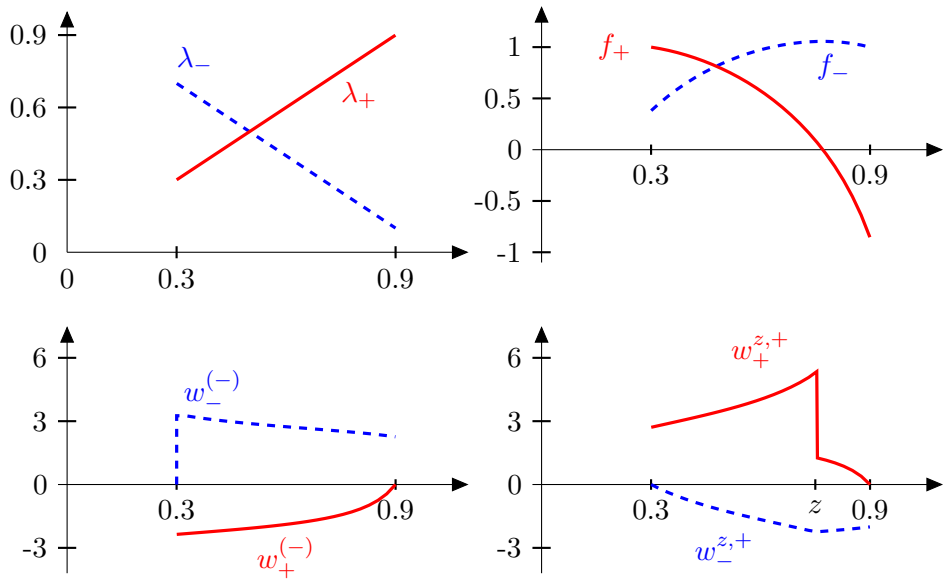

Figure~\ref{fig:numer} shows plots of weight functions as in
Theorem~\ref{T:spec} together with the solution of Luckock's equation, for one
explicit example of a Stigler-Luckock model. Theorem~\ref{T:spec} is closely
related to Luckock's result Theorem~\ref{T:Luckock}. Indeed, if a
Stigler-Luckock model has an invariant law that is concentrated on $\Si^{\rm
  fin}_{\rm ord}$, then the fact that the functions in (\ref{fpm}) are given
by the solution to Luckock's equation follows from Theorem~\ref{T:spec} and
the equilibrium equation $\E[GF(\Xc_t)]=0$.

Theorem~\ref{T:spec} is more powerful that Theorem~\ref{T:Luckock}, however,
since it gives an interpretation to the solution to Luckock's equation even if
such a solution is not valid. Also, we have fairly explicit expressions for
the weight functions $(w^{z,\pm}_-,w^{z,\pm}_+)$ (see Lemma~\ref{L:wint}
below), and their associated linear functions $F^{z,\pm}$ are useful also in a
non-stationary setting. In particular, we will prove Theorem~\ref{T:posrec} by
constructing a Lyapunov function from the functions $F^{I_-,-}$ and
$F^{I_+,+}$ (see formula (\ref{Ldef}) below).

We hope that the linear functions $F^{z,\pm}$ from Theorem~\ref{T:spec} will
also prove useful in future work aimed at resolving the open problems
mentioned in Section~\ref{S:open}. In Appendix~\ref{A:future}, we have
recorded some concrete ideas on how the functions $F^{z,\pm}$ could possibly
be used to attack the conjectures of Section~\ref{S:open}.

\subsection{Outline}

In Section~\ref{S:difeq}, we investigate two differential
equations: Luckock's equation (\ref{Luckock}) and a differential equation that
allows one to solve the weight functions $w_\pm$ in terms of the functions
$q_\pm$ from (\ref{qpm}). In particular, we prove Theorem~\ref{T:spec} in
Subsection~\ref{S:Luceq}, Proposition~\ref{P:Luck} in
Subsection~\ref{S:somex}, and Proposition~\ref{P:Phi} and
  Lemmas~\ref{L:JJ}, \ref{L:comp}, and \ref{L:unif} in
Subsection~\ref{S:restproof}.


After the preparatory work on the differential equations in
Section~\ref{S:difeq}, the analysis of the Markov chain, which is contained in
Section~\ref{S:Markov}, is actually quite short. In particular, we prove
Theorem~\ref{T:Luckock} in Subsection~\ref{S:station} and 
Theorem~\ref{T:posrec} in Subsection~\ref{S:posprf}.

The paper concludes with four appendices. In Appendix~\ref{A:staform}, we show
that the assumptions (A3) and (A4) can basically be made without loss of
generality.
Appendix~\ref{A:posrec} collects some facts from the general theory of Markov
chains needed to translate the properties of our Lyapunov function into
properties of the Markov chain. In Appendix~\ref{A:dis} we have collected
(without proof) some formulas for Stigler-Luckock models that take only
finitely many values, and that are analogues to our integral formulas for
continuous models but cannot easily be deduced from
them. Appendix~\ref{A:future} collects some concrete open problems with some
ideas on how to approach them.

\section{Analysis of the differential equations}\label{S:difeq}

\subsection{Lebesgue-Stieltjes integrals}\label{S:Stieltjes}

For any interval $J\sub[-\infty,\infty]$ that can be either closed, open, or
half open, with left and right boundaries $J_-<J_+$, we let $B(J)$ denote the
space of bounded measurable functions $f:J\to\R$. If a
function $f:J\to\R$ is of bounded variation, then the limits
\be
f(x-):=\lim_{y\up x}f(y)\quad(x\neq J_-)
\quand f(x+):=\lim_{y\down x}f(y)\quad(x\neq J_+)
\ee
exist for all $x\in J\beh\{J_-\}$ resp.\ $x\in J\beh\{J_+\}$. For such
functions, if $J_-\in J$, then we set $f(J_--):=f(J_-)$, and we define $f(J_++)$
similarly. We let $B_{\rm bv}(J)$ denote the space of functions
$f\in B(J)$ that are of bounded variation and satisfy
$f(x)\in\{f(x-),f(x+)\}$ for each $x\in J$, and we let
\be
B^\pm_{\rm bv}(J):=\big\{f\in B_{\rm bv}(J):f(x\pm)=f(x)\ \forall x\in J\big\}
\ee
denote the spaces of left $(-)$ and right $(+)$ continuous functions
$f:J\to\R$ of bounded variation. Each $f\in B_{\rm bv}(J)$ defines a finite
signed measure $\di f$ on $J$ through the formula
\be
\di f([x,y]):=f(y+)-f(x-)\qquad(x,y\in J,\ x\leq y).
\ee
The set of atoms of $\di f$ is the set $\Di_f:=\{x\in J:f(x-)\neq f(x+)\}$ of
points of discontinuity of $f$. For each finite signed measure $\rho$ on $J$
we can find functions $f\in B^-_{\rm bv}(J)$ and $g\in B^+_{\rm bv}(J)$ such
that $\di f=\rho=\di g$, and these functions are unique up to an additive
constant. We equip $B^\pm_{\rm bv}(J)$ with a topology such that $f_n\to f$ if
and only if $\di f_n$ converges weakly to $\di f$ and $f_n(x)\to f(x)$ for at
least one (and hence every) point $x\in J\beh\Di_f$. It is known
\cite[page~182]{Hoe77} that if $J$ is a closed interval, then $f_n\to f$ in
this topology if and only if:
\be\ba{rl}\label{Hoe77}
{\rm(i)}&\dis\sup_n\|\di f_n\|<\infty,
\mbox{ where $\|\,\cdot\,\|$ denotes the total variation norm},\\[5pt]
{\rm(ii)}&\dis\di f_n(J)\to\di f(J),\\[5pt]
{\rm(iii)}&\dis\int_J\big|f_n(x)-f(x)\big|\,\di x\to 0,
\mbox{ i.e., $f_n\to f$ in $L^1$ norm}\\
&\dis\mbox{w.r.t.\ to the Lebesgue measure.}
\ec



%

In line with earlier notation, we write $g\,\di f$ to
denote the measure $\di f$ weighted with a bounded measurable function $g$.
We will make use of the \emph{product rule} which says that
\be\label{prodrule}
\di(fg)=f\,\di g+g\,\di f
\qquad\big(f,g\in B_{\rm bv}(J),\ \Di_f\cap\Di_g=\emptyset\big),
\ee
and also of the \emph{chain rule} which tells us that if $f\in B_{\rm bv}(J)$
takes values in a compact interval $K$ and $F:K\to\R$ is continuously
differentiable, then
\be\label{chainrule}
\di(F\circ f)=(F'\circ f)\,\di f
\qquad\big(f\in B_{\rm bv}(J),\ \Di_f=\emptyset\big),
\ee
where $(F\circ f)(x):=F(f(x))$ denotes the composition of $F$ and $f$.
All our integrals will be of Lebesgue type, which coincides with the
Riemann-Stieltjes integral if both functions involved are of bounded variation
and do not share points of discontinuity.

If $g:[a,b]\to\R$ is a function of bounded variation and
$\psi:[a,b]\to\R$ is nondecreasing, then we write $\di
g\ll\di\psi$ if $\di g$ is absolutely continuous with respect to $\psi$, i.e.,
if for each $s\leq t$, $\psi(s-)=\psi(t+)$ implies $g(s-)=g(t+)$. We will
sometimes use the \emph{substitution of variables rule}, which says that
\be\label{substitute}
\int_a^b\!f\,\di g
=\int_{\psi(a)}^{\psi(b)}\!(f\circ\psi^{-1})\,\di(g\circ\psi^{-1})
\qquad\big(f\in B[a,b],\ g\in B_{\rm bv}[a,b]),
\ee
and which holds provided $\psi:[a,b]\to[-\infty,\infty]$ is a nondecreasing
function such that $\di g\ll\di\psi$, and
$\psi^{-1}:[\psi(a),\psi(b)]\to[a,b]$ is a right inverse of $\psi$.

As a general reference to these rules, we refer to
\cite[Section~6.2]{CB00}. In the substitution of variables rule, the condition
$\di g\ll\di\psi$ guarantees that $f\circ\psi^{-1}\circ\psi$ differs from $f$
only on a set of measure zero under $\di g$.

We will need one more result that we formulate as a lemma. The result holds in
any dimension but since we only need the two-dimensional case, for ease of
notation, we restrict to two dimensions.

\bl[Integrals along curves]\label{L:curve}
Let $D\sub\R^2$ be a closed, convex set that is the closure of its interior.
Let $F,g_1,g_2,f_1,f_2$ be continuous real functions on $D$ such that $f_1$
and $f_2$ are moreover Lipschitz. Assume that for all $x_1,x'_1,x_2,x'_2$ such
that $x_1\leq x'_1$, $x_2\leq x'_2$, and $(x_1,x_2),(x'_1,x_2),(x_1,x'_2)\in
D$,
\be\ba{r@{\,}c@{\,}ll}\label{pardif}
\dis F(x'_1,x_2)-F(x_1,x_2)
&=&\dis\int_{x_1}^{x'_1}\!g_1(\,\cdot\,,x_2)\,\di f_1(\,\cdot\,,x_2)\\[5pt]
\dis F(x_1,x'_2)-F(x_1,x_2)
&=&\dis\int_{x_2}^{x'_2}\!g_2(x_1,\,\cdot\,)\,\di f_2(x_1,\,\cdot\,).
\ec
Let $[t_-,t_+]$ be a closed interval and let $\ga:[t_-,t_+]\to D$ be a
continuous function of bounded variation. Then
\be\label{curve}
F\big(\ga(t_+)\big)-F\big(\ga(t_-)\big)
=\int_{t_-}^{t_+}\big\{(g_1\circ\ga)\,\di(f_1\circ\ga)
+(g_2\circ\ga)\,\di(f_2\circ\ga)\big\}.
\ee
\el
\bpro[Proof (sketch)]
Formula (\ref{pardif}) shows that (\ref{curve}) holds for any continuous
function $[t_-,t_+]\mapsto\big(\ga_1(t),\ga_2(t)\big)\in D$ of bounded
variation such that moreover either $\ga_1$ or $\ga_2$ is constant.  It
follows that (\ref{curve}) also holds for any finite concatenation of such
curves; call such curves simple. Then it is not hard to see that any
$\ga:[t_-,t_+]\to D$ that is continuous and of bounded variation can be
approximated by simple curves $\ga^{(n)}$ in such a way that
$\ga^{(n)}(t_-)\to\ga(t_-)$ and $\di\ga^{(n)}_i$
converges weakly to $\di\ga_i$ for $i=1,2$. In particular, this implies that
$\ga^{(n)}$ converges uniformly to $\ga$ so by the continuity of $g_i$
$(i=1,2)$, also $g_i\circ\ga^{(n)}$ converges uniformly to $g_i\circ\ga$.
In view of (\ref{Hoe77}), the Lipschitz continuity of $f_i$ $(i=1,2)$ moreover
implies that $\di(f_i\circ\ga^{(n)})$ converges weakly to $\di(f_i\circ\ga)$.
Using this and the continuity of $F$, taking the limit in (\ref{pardif}), which
holds for $\ga^{(n)}$, we obtain that the formula also holds for~$\ga$.
\epro

\subsection{The inverse problem}

The main result of the present subsection is the following theorem.

\bt[Inverse problem]
Assume\label{T:inverse} (A3) and (A6). Then, for each pair of functions $(g_-,g_+)$ with
$g_\pm\in B^\pm_{\rm bv}(\ov I)$, there exists a unique pair of functions
$(w^{(g_-,g_+)}_-,w^{(g_-,g_+)}_+)=(w_-,w_+)$ with $w_\pm\in B^\pm_{\rm
  bv}(\ov I)$ and $w_\pm(I_\pm)=0$, as well as a unique constant
$c(g_-,g_+)\in\R$, such that the linear functional $F^{(g_-,g_+)}=F$ from
(\ref{Flin}) satisfies
\be\label{spec}
GF^{(g_-,g_+)}(\Xc)=g_-\big(M_-(\Xc)\big)+g_+\big(M_+(\Xc)\big)-c(g_-,g_+).
\ee
\et

The proof of Theorem~\ref{T:inverse} will be split into a number of lemmas.

\bl[Differential equation]\label{L:wdif}
Assume (A3) and (A6), let $g_\pm\in B^\pm_{\rm bv}(\ov I)$, and let $w_\pm\in
B^\pm_{\rm bv}(\ov I)$ satisfy $w_\pm(I_\pm)=0$. Then the linear function
$F^{(w_-,w_+)}$ associated with $(w_-,w_+)$ satisfies (\ref{spec}) for some
$c(g_-,g_+)\in\R$ if and only if
\be\ba{rr@{\,}c@{\,}l}\label{wdif}
{\rm(i)}&\dis w_+\di\la_++\di(\la_+w_-)
&=&\dis-\di g_-,\\[5pt]
{\rm(ii)}&\dis w_-\di\la_-+\di(\la_-w_+)
&=&\dis-\di g_+.
\ec
\el
\bpro
Defining functions $q_\pm$ as in (\ref{qpm}), Lemma~\ref{L:GF}
tells us that (\ref{spec}) is satisfied for some $c(g_-,g_+)\in\R$ if and only
if there exist real constants $c_\pm$ such that $q_\pm=g_\pm+c_\pm$,
or equivalently, if there exist $c'_\pm\in\R$ such that
\be\ba{r@{\,}c@{\,}ll}
\dis g_-(x)&=&c'_-+\dis\int_{[x,I_+)}\big\{w_+\di\la_++\di(w_-\la_+)\big\}
\qquad\big(x\in[I_-,I_+)\big),\\[5pt]
\dis g_+(x)&=&\dis c'_+-\int_{(I_-,x]}\big\{w_-\di\la_-+\di(w_+\la_-)\big\}
\qquad\big(x\in(I_-,I_+]\big),
\ec
which is equivalent to (\ref{wdif}).
\epro

We can integrate the differential equation (\ref{wdif}) explicitly. Let
$h_\pm\in B^\pm_{\rm bv}(\ov I)$ be any pair of functions such that
\be\label{hprop}
\di h_\pm=-\la_\pm\di g_\pm,
\ee
i.e., $h_-(x)=c_-+\int_{[x,I_+)}\la_-\di g_-$ and similarly for $h_+$,
where $c_\pm$ are some fixed, but otherwise arbitrary constants.

\bl[Integrated equation]
Assume\label{L:Lucinh} (A3) and (A6), let $g_\pm\in B^\pm_{\rm bv}(\ov I)$ be given and let
$h_\pm\in B^\pm_{\rm bv}(\ov I)$ be as in (\ref{hprop}). Then a pair of
functions $w_\pm\in B^\pm_{\rm bv}(\ov I)$ satisfies
(\ref{wdif}) if and only if there exists a constant $\kappa\in\R$ such that
\be\ba{rr@{\,}c@{\,}l}\label{sol}
{\rm(i)}&\dis\di w_-
&=&\dis\frac{\kappa+h_+}{\la_-}\,\di\Big(\frac{1}{\la_+}\Big)
+\frac{1}{\la_-}\,\di\Big(\frac{h_-}{\la_+}\Big),\\[7pt]
{\rm(ii)}&\dis\di w_+
&=&\dis\frac{\kappa+h_-}{\la_+}\,\di\Big(\frac{1}{\la_-}\Big)
+\frac{1}{\la_+}\,\di\Big(\frac{h_+}{\la_-}\Big),\\[7pt]
{\rm(iii)}&\dis w_-+w_+
&=&\dis\frac{\kappa+h_-+h_+}{\la_-\la_+}.
\ec
Moreover, given (\ref{sol})~(iii), the equations (\ref{sol})~(i)
and (ii) imply each other.
\el
\bpro
Multiplying the equations (\ref{wdif})~(i) and (ii) by $\la_-$ and $\la_+$,
respectively, and then adding both equations, using the product rule
(\ref{prodrule}) and (\ref{hprop}), we obtain
\be
\di\big(\la_-(\la_+w_-)\big)+\di\big(\la_+(\la_-w_+)\big)
=\di h_-+\di h_+,
\ee
which shows that there exists a constant
$\kappa\in\R$ such that (\ref{sol})~(iii) holds.
Given (\ref{sol})~(iii), we can rewrite
(\ref{wdif})~(i) as
\be
\di(\la_+w_-)=-w_+\di\la_++\frac{\di h_-}{\la_-}
=\Big(w_--\frac{\kappa+h_-+h_+}{\la_-\la_+}\Big)\di\la_+
+\frac{\di h_-}{\la_-}.
\ee
Dividing by $\la_+$ and reordering terms, this says that
\be
\frac{\di(\la_+w_-)-w_-\di\la_+}{\la_+}
=-\frac{(\kappa+h_-+h_+)\,\di\la_+}{\la_-\la_+^2}
+\frac{\di h_-}{\la_-\la_+},
\ee
which using the product and chain rules
(\ref{prodrule})--(\ref{chainrule}) can be rewritten as
(\ref{sol})~(i). In a similar way, we see that given (\ref{sol})~(iii),
(\ref{Luckock})~(ii) is equivalent to (\ref{sol})~(ii). Differentiating
(\ref{sol})~(iii), using the product rule, we obtain
\bc
\dis\di w_-+\di w_+
&=&\dis\frac{\kappa}{\la_-}\,\di\Big(\frac{1}{\la_+}\Big)
+\frac{\kappa}{\la_+}\,\di\Big(\frac{1}{\la_-}\Big)
+\frac{1}{\la_-}\,\di\Big(\frac{h_-}{\la_+}\Big)\\[5pt]
&&\dis+\frac{h_-}{\la_+}\,\di\Big(\frac{1}{\la_-}\Big)
+\frac{h_+}{\la_-}\,\di\Big(\frac{1}{\la_+}\Big)
+\frac{1}{\la_+}\,\di\Big(\frac{h_+}{\la_-}\Big),
\ec
which is the same as we would obtain adding the equations (\ref{sol})~(i)
and (\ref{sol})~(ii). We conclude that (\ref{sol})~(i)
and (ii) are equivalent given (iii).
\epro

For later use, assuming (A3) and (A6), we define a constant $\Ga$ by
\be\label{Gamma}
\Ga:=\frac{1}{\la_-(I_+)\la_+(I_+)}
-\int_{I_-}^{I_+}\frac{1}{\la_-}\di\Big(\frac{1}{\la_+}\Big)
=\frac{1}{\la_-(I_-)\la_+(I_-)}
+\int_{I_-}^{I_+}\frac{1}{\la_+}\di\Big(\frac{1}{\la_-}\Big),
\ee
where the equality of both formulas follows from the product rule
(\ref{prodrule}) applied to the functions $1/\la_-$ and $1/\la_+$. Note that
$\Ga>0$ since $\di(1/\la_-)$ is nonnegative while $\la_\pm$ are strictly
positive by (A6).

\bl[Existence and uniqueness]\label{L:exuni}
Assume (A3) and (A6). Then, for each pair of functions $g_\pm\in B^\pm_{\rm
  bv}(\ov I)$, there exist unique functions $w_\pm\in B^\pm_{\rm bv}(\ov I)$
that solve the differential equation (\ref{wdif}) together with
the boundary conditions $w_\pm(I_\pm)=0$. These functions are given by
\be\ba{rr@{\,}c@{\,}l}\label{wsol}
{\rm(i)}&\dis w_-(x)&=&\dis\int_{[I_-,x)}
\Big\{\frac{\kappa+h_+}{\la_-}\,\di\Big(\frac{1}{\la_+}\Big)
+\frac{1}{\la_-}\,\di\Big(\frac{h_-}{\la_+}\Big)\Big\},\\[5pt]
{\rm(ii)}&\dis w_+(x)&=&\dis-\int_{(x,I_+]}
\Big\{\frac{\kappa+h_-}{\la_+}\,\di\Big(\frac{1}{\la_-}\Big)
+\frac{1}{\la_+}\,\di\Big(\frac{h_+}{\la_-}\Big)\Big\},
\ec
where $h_\pm$ are as in (\ref{hprop}) and
\be\label{kappah}
\kappa
:=\Ga^{-1}\Big[
\int_{\ov I}\Big\{\frac{h_+}{\la_-}\,\di\Big(\frac{1}{\la_+}\Big)
+\frac{1}{\la_-}\,\di\Big(\frac{h_-}{\la_+}\Big)\Big\}
-\frac{h_-(I_+)+h_+(I_+)}{\la_-(I_+)\la_+(I_+)}\Big],
\ee
with $\Ga>0$ the constant from (\ref{Gamma}).
\el
\bpro
By Lemma~\ref{L:Lucinh}, $w_\pm$ solve the difference equation
(\ref{wdif}) together with the left boundary condition $w_-(I_-)=0$
if and only if there exists a $\kappa\in\R$ such that (\ref{wsol})~(i)
and (\ref{sol})~(iii) hold. In view of the latter equation, $w_\pm$
also solves the right boundary condition $w_+(I_+)=0$ if and only if
\be
w_-(I_+)+0=\frac{\kappa+h_-+h_+}{\la_-\la_+}(I_+).
\ee
In view of (\ref{wsol})~(i), this says that
\be
\int_{[I_-,I_+)}\Big\{\frac{\kappa+h_+}{\la_-}\,\di\Big(\frac{1}{\la_+}\Big)
+\frac{1}{\la_-}\,\di\Big(\frac{h_-}{\la_+}\Big)\Big\}
=\frac{\kappa+h_-(I_+)+h_+(I_+)}{\la_-(I_+)\la_+(I_+)},
\ee
or equivalently (note that since $h_-$ is left-continuous, it has no jump at
$I_+$)
\be\ba{l}
\dis\int_{[I_-,I_+]}\Big\{\frac{h_+}{\la_-}\,\di\Big(\frac{1}{\la_+}\Big)
+\frac{1}{\la_-}\,\di\Big(\frac{h_-}{\la_+}\Big)\Big\}
-\frac{h_-(I_+)+h_+(I_+)}{\la_-(I_+)\la_+(I_+)}\\[5pt]
\dis\qquad\qquad=\kappa\Big\{\frac{1}{\la_-(I_+)\la_+(I_+)}
-\int_{I_-}^{I_+}\frac{1}{\la_-}\,\di\Big(\frac{1}{\la_+}\Big)\Big\},
\ec
which by the fact that the constant $\Ga$ from (\ref{Gamma}) is nonzero
is equivalent to (\ref{kappah}).
\epro

\bpro[Proof of Theorem~\ref{T:inverse}]
Immediate from Lemmas~\ref{L:wdif} and \ref{L:exuni}.
\epro

\subsection{Luckock's equation}\label{S:Luceq}

In the present subsection we prove Theorem~\ref{T:spec}. We start by proving
that Luckock's equation has a unique solution. By definition, a solution to
\emph{Luckcock's equation} is a pair of functions $(f_-,f_+)$ such that
$f_\pm\in B^\mp_{\rm bv}(\ov I)$ and (\ref{Luckock}) holds.

\bl[Luckock's equation]\label{L:Lucsol}
Assume (A3) and (A6). Then Luckock's equation has a unique solution
$(f_-,f_+)$, which is given by
\be\ba{rr@{\,}c@{\,}l}\label{Lucsol}
{\rm(i)}&\dis\Big(\frac{f_+}{\la_+}\Big)(x)
&=&\dis\frac{1}{\la_+(I_-)}
+\kappa\int_{I_-}^x\frac{1}{\la_-}\,\di\Big(\frac{1}{\la_+}\Big),\\[5pt]
{\rm(ii)}&\dis\Big(\frac{f_-}{\la_-}\Big)(x)
&=&\dis\frac{1}{\la_-(I_+)}
-\kappa\int_x^{I_+}\frac{1}{\la_+}\,\di\Big(\frac{1}{\la_-}\Big),
\ec
where $\kappa$ is given by
\be\label{Luckappa}
\kappa=\kappa_{\rm L}:=\Ga^{-1}\Big(\frac{1}{\la_-(I_+)}+\frac{1}{\la_+(I_-)}\Big),
\ee
and $\Ga>0$ is the constant from (\ref{Gamma}).
\el
\bpro
Setting $v_+:=f_-/\la_-$ and $v_-:=f_+/\la_+$ and dividing the equations
(\ref{Luckock})~(i) and (ii) by $\la_-$ and $\la_+$, respectively, we see
that these equations are equivalent to $v_\pm\di\la_\pm=-\di(\la_\pm v_\mp)$,
which is equation (\ref{wdif}) with $w_\pm=v_\pm$ and $g_\pm=0$. Now
Lemma~\ref{L:Lucinh} tells us that $(f_-,f_+)$ solves
(\ref{Luckock})~(i) and (ii) if and only if there exists a constant
$\kappa\in\R$ such that
\be\ba{rr@{\,}c@{\,}l}\label{Lucint}
{\rm(i)}&\dis\di\Big(\frac{f_+}{\la_+}\Big)
&=&\dis\frac{\kappa}{\la_-}\,\di\Big(\frac{1}{\la_+}\Big),\\[7pt]
{\rm(ii)}&\dis\di\Big(\frac{f_-}{\la_-}\Big)
&=&\dis\frac{\kappa}{\la_+}\,\di\Big(\frac{1}{\la_-}\Big),\\[7pt]
{\rm(iii)}&\dis\frac{f_+}{\la_+}+\frac{f_-}{\la_-}
&=&\dis\frac{\kappa}{\la_-\la_+}.
\ec
Moreover, of these equations, the first two are equivalent given the third one.

It follows that $(f_-,f_+)$ solves (\ref{Luckock})~(i) and (ii) together
with the left boundary condition $f_+(I_-)=1$ if and only if
(\ref{Lucsol})~(i) and (\ref{Lucint})~(iii) hold. In view of the latter
equation, the right boundary condition $f_-(I_+)=1$ is satisfied if and only
if
\be
\frac{f_+}{\la_+}(I_+)+\frac{1}{\la_-(I_+)}
=\frac{\kappa}{\la_-\la_+}(I_+).
\ee
By (\ref{Lucsol})~(i), this says that
\be
\frac{1}{\la_+(I_-)}
+\kappa\int_{I_-}^{I_+}\frac{1}{\la_-}\,\di\Big(\frac{1}{\la_+}\Big)
+\frac{1}{\la_-(I_+)}=\frac{\kappa}{\la_-(I_+)\la_+(I_+)},
\ee
which is equivalent to (\ref{Luckappa}).
\epro

\bpro[Proof of Theorem~\ref{T:spec}]
Let $\Gi$ be the space of all pairs $(g_-,g_+)$ with $g_\pm\in B^\pm_{\rm
  bv}(\ov I)$ and set $\Wi:=\{(w_-,w_+)\in\Gi:w_\pm(I_\pm)=0\}$. We equip the
spaces $B^\pm_{\rm bv}(\ov I)$ with a topology as in
Section~\ref{S:Stieltjes}, $\Gi$ with the product topology, and $\Wi$ with the
induced topology. For any interval $J$, we let $\Mi(J)$ denote the space of 
finite signed measures on $J$, equipped with the topology of weak convergence,
and we let $\Ri:=\Mi[I_-,I_+)\times\Mi(I_-,I_+]$, equipped with the product
topology.

Let $\psi:\Wi\to\Gi$ be the linear function that maps a pair $(w_-,w_+)\in\Wi$
into the pair $(q_-,q_+)\in\Gi$ defined in (\ref{qpm}) and let $D:\Gi\to\Ri$
be the map
\be
D(g_-,g_+):=\big(\di g_-,\di g_+).
\ee
Setting $\phi:=D\circ\psi$, we see that
\be\label{phidef}
\phi(w_-,w_+)=-\big(w_+\di\la_++\di(\la_+w_-),
w_-\di\la_-+\di(\la_-w_+)\big),
\ee
so Lemma~\ref{L:exuni} tells us that $\phi:\Wi\to\Ri$ is a bijection.

We claim that the maps $\psi$, $D$, $\phi$, and $\phi^{-1}$ are continuous
with respect to the topologies on $\Wi$, $\Gi$, and $\Ri$. The continuity of
$D$ is immediate from the definition of the topologies on $\Gi$ and $\Ri$ and
the continuity of $\psi$ follows from (\ref{qpm}). The continuity of $\phi$ is
easily derived from (\ref{phidef}), while the continuity of $\phi^{-1}$
follows from the explicit formulas in Lemma~\ref{L:exuni} and the continuity
of the functions $h_\pm$ from (\ref{hprop}) as a function of $g_\pm$, for a
given choice of the boundary conditions.


Let $\psi(\Wi)$ denote the image of $\Wi$ under $\psi$ and define
$\pi:\Gi\to\psi(\Wi)$ by $\pi:=\psi\circ\phi^{-1}\circ D$. Since
$\pi\circ\psi=\psi\circ\phi^{-1}\circ(D\circ\psi)=\psi$, we see that $\pi$ is
the identity on $\psi(\Wi)$. Since $D\circ\pi=(D\circ\psi)\circ\phi^{-1}\circ
D=D$, we see that $\pi(g)=\pi(g')$ if and only if $D(g)=D(g')$. These facts
imply that for each $g\in\Gi$, there exists a unique $q\in\psi(\Wi)$, namely
$q=\pi(g)$, such that $D(g)=D(q)$, i.e., for every $g\in\Gi$ there exists a
unique $q\in\psi(\Wi)$ and unique constants $c_\pm(g_-,g_+)\in\R$ such that
\be\label{qg}
(g_-,g_+)=\big(q_-+c_-(g_-,g_+),q_++c_+(g_-,g_+)\big).
\ee
Since $\psi$, $\phi^{-1}$ and $D$ are continuous, so is $\pi$ and hence also
the maps $c_\pm:\Gi\to\R$ are continuous. In fact, they are the unique
continuous linear forms on $\Gi$ such that
\be\ba{rr@{\,}c@{\,}lr@{\,}c@{\,}lr@{\,}c@{\,}ll}
{\rm(i)}&c_-(1,0)&=&1,\ &c_-(0,1)&=&0,\ 
&c_-\big(\psi(w_-,w_+)\big)&=&0\ &\forall(w_-,w_+)\in\Wi,\\[5pt]
{\rm(ii)}&c_+(1,0)&=&0,\ &c_+(0,1)&=&1,\ 
&c_+\big(\psi(w_-,w_+)\big)&=&0\ &\forall(w_-,w_+)\in\Wi.
\ec
The map $(g_-,g_+)\mapsto c(g_-,g_+)$ from Theorem~\ref{T:inverse}
is given by $c=c_-+c_+$, i.e., $c$ is the unique
continuous linear form on $\Gi$ such that
\be\label{clin}
c(1,0)=1,\quad c(0,1)=1,\quad
c\big(\psi(w_-,w_+)\big)=0\quad\forall(w_-,w_+)\in\Wi.
\ee
Let $(f_-,f_+)$ be the unique solution to Luckock's equation, and observe from
(\ref{Lucsol}) that $f_\pm$ are continuous on $\ov I$. We claim that
\bc\label{cgg}
\dis c(g_-,g_+)
&=&\dis g_-(I_+)f_-(I_+)-\int_{\ov I}f_-\di g_-
+g_+(I_-)f_+(I_-)+\int_{\ov I}f_+\di g_+.
\ec
Clearly, (\ref{cgg}) defines a continuous linear form on $\Gi$.
We will show that this linear form satisfies (\ref{clin}).
The boundary conditions (\ref{Luckock})~(iii) imply
that $c(1,0)=1=c(0,1)$. Recall that for $(w_-,w_+)\in\Wi$,
$\psi(w_-,w_+)=(q_-,q_+)$ is defined as in (\ref{qpm}). Then
\bc\label{cpsi}
\dis c\big(\psi(w_-,w_+)\big)
&=&\dis-(w_-\la_+)(I_+)f_-(I_+)
+\int_{\ov I}f_-\big\{w_+\di\la_++\di(w_-\la_+)\big\}\\[5pt]
&&\dis-(w_+\la_-)(I_-)f_+(I_-)
-\int_{\ov I}f_+\big\{w_-\di\la_-+\di(w_+\la_-)\big\}.
\ec
By partial integration, using the continuity of $f_\pm$ and $\la_\pm$,
as well as the boundary condition $w_-(I_-)=0$, we have
\be
-(w_-\la_+)(I_+)f_-(I_+)+\int_{\ov I}f_-\di(w_-\la_+)
=
-\int_{\ov I}(w_-\la_+)\di f_-,
\ee
Inserting this into
the first line of (\ref{cpsi}) and treating the second line similarly, we find
that $c\big(\psi(w_-,w_+)\big)$ equals
\be\ba{l}\label{lastep}
\dis\int_{\ov I}\big\{f_-w_+\di\la_+-(w_-\la_+)\di f_-\big\}
+\int_{\ov I}\big\{-f_+w_-\di\la_-+(w_+\la_-)\di f_+\big\}\\[5pt]
\dis\quad=\int_{\ov I}w_+\big\{f_-\di\la_++\la_-\di f_+\big\}
-\int_{\ov I}w_-\big\{f_+\di\la_-+\la_+\di f_-\big\}=0,
\ec
where we have used (\ref{Luckock})~(i) and (ii) in the last step. This
completes the proof of (\ref{cgg}).

In particular, formula (\ref{cgg}) shows that
\be
c(1_{[I_-,z]},0)=f_-(z)
\quand
c(0,1_{[z,I_+]})=f_+(z)\qquad(z\in\ov I),
\ee
which together with Theorem~\ref{T:inverse} implies Theorem~\ref{T:spec}.
\epro

\subsection{Some explicit formulas and conditions}\label{S:somex}

In the present section, we prove Proposition~\ref{P:Luck} as well as two
lemmas (Lemmas~\ref{L:wint} and \ref{L:wextr} below) giving explicit formulas
for the weight functions of Theorem~\ref{T:spec}.\med



\bpro[Proof of Proposition~\ref{P:Luck}]
The fact that under the conditions (A3) and (A6), Luckock's equation has a
unique solution has already been proved in Lemma~\ref{L:Lucsol}.

In order to prove (\ref{valid}), it suffices to prove part~(i); the other part
then follows by symmetry. Let
$\La_{+-}:=\int_{I_-}^{I_+}\frac{1}{\la_+}\di\big(\frac{1}{\la_-}\big)$.
Then (\ref{Lucsol})~(ii) says that
\be
f_-(I_-)=\la_-(I_-)\Big\{\frac{1}{\la_-(I_+)}-\kappa_L\La_{+-}\Big\}.
\ee
Filling in the definition of $\kappa_L$ in (\ref{Luckappa}), we see that
$f_-(I_-)>0$ if and only if
\be
\frac{1}{\la_-(I_+)}>
\Ga^{-1}\Big(\frac{1}{\la_-(I_+)}+\frac{1}{\la_+(I_-)}\Big)\La_{+-}.
\ee
Using also formula (\ref{Gamma}) and the fact that $\Ga>0$, this can be
rewritten as 
\be
\Big(\frac{1}{\la_-(I_-)\la_+(I_-)}+\La_{+-}\Big)\frac{1}{\la_-(I_+)}
>\Big(\frac{1}{\la_-(I_+)}+\frac{1}{\la_+(I_-)}\Big)\La_{+-},
\ee
which can be simplified to (\ref{valid})~(i). The same argument also works
with all inequality signs reversed.

To complete the proof, we need to show (\ref{volume}). Consider the weight
functions
\be\label{wmass}
w_-:=-1_{(I_-,I_+]}\quand w_+:=1_{[I_-,I_+)},
\ee
which correspond through (\ref{Flin}) to the linear function $F(\Xc)=\Xc(I)$.
For these weight functions, the functions $q_\pm$ from (\ref{qpm}) are given
by
\bc
\dis q_-(x)&=&\dis
\la_+(I_+)-\la_+(I_-)1_{\{I_-\}}(x),\\[5pt]
\dis q_+(x)&=&\dis
-\la_-(I_-)+\la_-(I_+)1_{\{I_+\}}(x),
\ec
so Lemma~\ref{L:GF} tells us that
\be
GF(\Xc)=-\la_+(I_-)1_{\{M_-(\Xc)=I_-\}}+\la_-(I_+)1_{\{M_+(\Xc)=I_+\}}
+\la_+(I_+)-\la_-(I_-).
\ee
By Theorem~\ref{T:inverse}, the weight functions $w_\pm$ are in fact uniquely
characterized by the requirement that
\be\label{plusc}
GF(\Xc)=-\la_+(I_-)1_{\{M_-(\Xc)=I_-\}}+\la_-(I_+)1_{\{M_+(\Xc)=I_+\}}+c
\ee
for some $c\in\R$. Defining weight functions $\ti w_\pm$ by
\be\label{tiwmass}
\ti w_\pm:=-\la_+(I_-)w^{I_-,-}_\pm+\la_-(I_+)w^{I_+,+}_\pm
\ee
and denoting the corresponding linear function by $\ti F$, we see from 
Theorem~\ref{T:spec} that
\bc
\dis G\ti F(\Xc)
&=&\dis-\la_+(I_-)1_{\{M_-(\Xc)=I_-\}}+\la_-(I_+)1_{\{M_+(\Xc)=I_+\}}\\[5pt]
&&\dis+\la_+(I_-)f_-(I_-)-\la_-(I_+)f_+(I_+).
\ec
We conclude from this that $w_\pm=\ti w_\pm$ and the constant from
(\ref{plusc}) is given by
\be
\la_+(I_+)-\la_-(I_-)=c=\la_+(I_-)f_-(I_-)-\la_-(I_+)f_+(I_+),
\ee
which proves (\ref{volume}).
\epro

We next set out to derive explicit formulas for the weight functions
$(w^{z,\pm}_-,w^{z,\pm}_+)$ from Theorem~\ref{T:spec}. To state the result, we
define functions $\Vv_{\pm,\mp}:\ov I\to\R$ by
\bc\label{Vpm}
\Vv_{-+}(x)&:=&\dis\Ga^{-1}\Big\{\frac{1}{\la_-(I_+)\la_+(I_+)}
-\int_x^{I_+}\frac{1}{\la_-}\di\Big(\frac{1}{\la_+}\Big)\Big\},\\[5pt]
\Vv_{+-}(x)&:=&\dis\Ga^{-1}\Big\{\frac{1}{\la_-(I_-)\la_+(I_-)}
+\int_{I_-}^x\frac{1}{\la_+}\di\Big(\frac{1}{\la_-}\Big)\Big\}.
\ec
In view of (\ref{Gamma}), we observe that
$\Vv_{-+}(I_-)=1=\Vv_{+-}(I_+)$. Moreover, $\Vv_{-+}$ is nonincreasing with
$\Vv_{-+}(I_+)>0$ while $\Vv_{+-}$ is nondecreasing with $\Vv_{-+}(I_-)>0$.
By partial integration, our formulas for $\Vv_{-+}$ and $\Vv_{+-}$ can be
rewritten as
\bc\label{Vpm2}
\Vv_{-+}(x)&:=&\dis\Ga^{-1}\Big\{\frac{1}{\la_-(x)\la_+(x)}
+\int_x^{I_+}\frac{1}{\la_+}\di\Big(\frac{1}{\la_-}\Big)\Big\},\\[5pt]
\Vv_{+-}(x)&:=&\dis\Ga^{-1}\Big\{\frac{1}{\la_-(x)\la_+(x)}
-\int_{I_-}^x\frac{1}{\la_-}\di\Big(\frac{1}{\la_+}\Big)\Big\}.
\ec
Combining this with our previous formulas and (\ref{Gamma}), we see that
\be\label{VV}
\Vv_{-+}(x)+\Vv_{+-}(x)
=\frac{\Ga^{-1}}{\la_-(x)\la_+(x)}+1.
\ee

\bl[Formulas for special weight functions]
The\label{L:wint} weight functions from Theorem~\ref{T:spec} are given by
\be\ba{rr@{\,}c@{\,}l}\label{wint}
{\rm(i)}&\dis w^{z,-}_-(x)&=&\dis\la_-(z)\Ga\big(\Vv_{+-}(z)-1_{\{x\leq z\}}\big)
\big(\Vv_{-+}(x)-1_{\{x\leq z\}}\big)\\[5pt]
{\rm(ii)}&\dis w^{z,-}_+(x)
&=&\dis\la_-(z)\Ga\big[\Vv_{+-}(x\vee z)-1\big]\Vv_{+-}(x\wedge z)\\[5pt]
{\rm(iii)}&\dis w^{z,+}_-(x)
&=&\dis\la_+(z)\Ga\big[\Vv_{-+}(x\wedge z)-1\big]\Vv_{-+}(x\vee z)\\[5pt]
{\rm(iv)}&\dis w^{z,+}_+(x)&=&\dis\la_+(z)\Ga\big(\Vv_{-+}(z)-1_{\{x\geq z\}}\big)
\big(\Vv_{+-}(x)-1_{\{x\geq z\}}\big).
\ec
\el
\bpro
We start with formula (\ref{wint})~(ii). Since $w^{I_+,-}_+=0$ which agrees
with the right-hand side of (\ref{wint})~(ii), we assume from now on without
loss of generality that $z\in[I_-,I_+)$. We apply Lemma~\ref{L:exuni} with
$g_-=1_{[I_-,z]}$ and $g_+=0$. For the functions $h_\pm$ from (\ref{hprop}) we
choose the boundary conditions $h_-(I_+)=0=h_+(I_-)$, which means that
\be\label{hz}
h_-(x)=\int_{[x,I_+)}\la_-\di 1_{[I_-,z]}=-\la_-(z)1_{[I_-,z]}(x)
\quand
h_+=0.
\ee
Since $h_+=0$ and $h_-(I_+)=0$, formulas (\ref{wsol})~(ii) and (\ref{kappah})
now simplify to
\be
w_+(x)=-\int_{(x,I_+]}
\frac{\kappa+h_-}{\la_+}\,\di\Big(\frac{1}{\la_-}\Big)
\quad\mbox{with}\quad
\kappa=\Ga^{-1}\int_{\ov I}\frac{1}{\la_-}\,\di\Big(\frac{h_-}{\la_+}\Big).
\ee
Here, by (\ref{hz}),
\bc
\dis\int_{\ov I}\frac{1}{\la_-}\,\di\Big(\frac{h_-}{\la_+}\Big)
&=&\dis\Big[\frac{h_-}{\la_-\la_+}(I_+)-\frac{h_-}{\la_-\la_+}(I_-)\Big]
-\int_{\ov I}\frac{h_-}{\la_+}\,\di\Big(\frac{1}{\la_-}\Big)\\[5pt]
&=&\dis\la_-(z)\Big\{\frac{1}{\la_-(I_-)\la_+(I_-)}
-\int_{[I_-,z]}\frac{1}{\la_+}\,\di\Big(\frac{1}{\la_-}\Big)\Big\},
\ec
which shows that
\be\label{klV}
\kappa=\la_-(z)\Vv_{+-}(z).
\ee
Using the fact that $\Vv_{+-}(I_+)=1$, it follows that
\bc
\dis w_+(x)
&=&\dis-\la_-(z)\Vv_{+-}(z)
\int_{(x,I_+]}\frac{1}{\la_+}\,\di\Big(\frac{1}{\la_-}\Big)
+\la_-(z)1_{\{x<z\}}\int_{(x,z]}\frac{1}{\la_+}\,\di\Big(\frac{1}{\la_-}\Big)\\[5pt]
&=&\dis-\la_-(z)\Ga\Big\{
\Vv_{+-}(z)\big[1-\Vv_{+-}(x)\big]-1_{\{x<z\}}\big[\Vv_{+-}(z)-\Vv_{+-}(x)\big]\Big\},
\ec
which can be rewritten as (\ref{wint})~(ii).

We next prove (\ref{wint})~(i). By (\ref{sol})~(iii), (\ref{hz}),
(\ref{klV}), and (\ref{wint})~(ii),
\bc
\dis w_+(x)&=&\dis\frac{\kappa+h_-+h_+}{\la_-\la_+}(x)-w_-(x)\\[5pt]
&=&\dis\la_-(z)\frac{\Vv_{+-}(z)-1_{[I_-,z]}}{\la_-\la_+}
-\la_-(z)\Ga\big(\Vv_{+-}(x\vee z)-1\big)\Vv_{+-}(x\wedge z).
\ec
For $x\leq z$, using (\ref{VV}), this yields
\bc
\dis w_+(x)
&=&\dis\la_-(z)\Ga\Big\{\Ga^{-1}\frac{\Vv_{+-}(z)-1}{\la_-\la_+}
-\big(\Vv_{+-}(z)-1\big)\Vv_{+-}(x)\Big\}\\[5pt]
&=&\dis\la_-(z)\Ga\big(\Vv_{+-}(z)-1\big)
\Big\{\frac{\Ga^{-1}}{\la_-\la_+}-\Vv_{+-}(x)\Big\}\\[5pt]
&=&\dis\la_-(z)\Ga\big(\Vv_{+-}(z)-1\big)
\big(\Vv_{-+}(x)-1\big),
\ec
while for $x>z$, again with the help of (\ref{VV}), we obtain
\bc
\dis w_+(x)
&=&\dis\la_-(z)\Ga\Big\{\Ga^{-1}\frac{\Vv_{+-}(z)}{\la_-\la_+}
-\big(\Vv_{+-}(x)-1\big)\Vv_{+-}(z)\Big\}\\[5pt]
&=&\dis\la_-(z)\Ga \Vv_{+-}(z)\Big\{\frac{\Ga^{-1}}{\la_-\la_+}
-\Vv_{+-}(x)+1\Big\}\\[5pt]
&=&\dis\la_-(z)\Ga \Vv_{+-}(z)\Vv_{-+}(x).
\ec
Combining the previous two formulas, we arrive at (\ref{wint})~(i).

Formulas (\ref{wint})~(iii) and (\ref{wint})~(iv) can be proved in exactly the
same way. Alternatively, they can be derived from (\ref{wint})~(ii) and
(\ref{wint})~(i) using the symmetry between buy and sell orders.
\epro

Assume (A3) and (A6) and for $z\in\ov I$, let $F^{z,\pm}$ be linear
functionals defined in terms of weight functions
$(w^{z,\pm}_-,w^{z,\pm}_+)$ as in Theorem~\ref{T:spec}.
We will in particular be interested in the case $z=I_\pm$ and introduce the
shorthands
\be\label{wbp}
w^{(\pm)}_-:=w^{I_\pm,\pm}_-\quad w^{(\pm)}_+:=w^{I_\pm,\pm}_+,
\quand F^{(\pm)}:=F^{I_\pm,\pm}.
\ee
We will prove Theorem~\ref{T:posrec} by constructing a Lyapunov function from
$F^{(-)}$ and $F^{(+)}$, see formula (\ref{Ldef}) and Proposition~\ref{P:Lyap}
below. The next lemma prepares for the proof of Proposition~\ref{P:Lyap}.

\bl[Extremal weight functions]
Assume\label{L:wextr} (A3) and (A6), let $(f_-,f_+)$ denote the solution to Luckock's
equation, and let $(w^{(\pm)}_-,w^{(\pm)}_+)$ be defined as in
(\ref{wbp}). Then for $x\in I$, one has
\be\ba{rr@{\,}c@{\,}l@{\quad}rr@{\,}c@{\,}l}\label{wextr}
{\rm(i)}&\dis w^{(-)}_-(x)&=&\dis\frac{\Vv_{-+}(x)}{\la_+(I_-)},
&{\rm(ii)}&\dis w^{(-)}_+(x)
&=&\dis-\frac{1-\Vv_{+-}(x)}{\la_+(I_-)},\\[10pt]
{\rm(iii)}&\dis w^{(+)}_-(x)&=&\dis-\frac{1-\Vv_{-+}(x)}{\la_-(I_+)},
&{\rm(iv)}&\dis w^{(+)}_+(x)&=&\dis\frac{\Vv_{+-}(x)}{\la_-(I_+)}.
\ec
Moreover,
\be\ba{rrcl}\label{wsum}
{\rm(i)}&\dis f_+(I_+)>0&\desd
&\dis\inf_{x\in I}\big[w^{(-)}_-(x)+w^{(+)}_-(x)\big]>0,\\[5pt]
{\rm(ii)}&\dis f_-(I_-)>0&\desd
&\dis\inf_{x\in I}\big[w^{(-)}_+(x)+w^{(+)}_+(x)\big]>0.
\ec
Both formulas also hold with the inequality signs reversed.
\el
\bpro
We only prove (\ref{wextr})~(i) and (ii) and (\ref{wsum})~(i); the proof of
the other formulas follows from the symmetry between buy and sell orders.
By Lemma~\ref{L:wint} and the facts that
\be
\Vv_{-+}(I_+)=\frac{\Ga^{-1}}{\la_-(I_+)\la_+(I_+)}
\quand
\Vv_{+-}(I_-)=\frac{\Ga^{-1}}{\la_-(I_-)\la_+(I_-)},
\ee
we have
\be\ba{rr@{\,}c@{\,}l}
{\rm(i)}&\dis w^{(-)}_-(x)&=&\dis\la_-(I_-)\Ga
\Big[\frac{\Ga^{-1}}{\la_-(I_-)\la_+(I_-)}-1_{\{x=I_-\}}\Big]
\big(\Vv_{-+}(x)-1_{\{x=I_-\}}\big)\\[5pt]
{\rm(ii)}&\dis w^{(-)}_+(x)
&=&\dis\la_-(I_-)\Ga\big[\Vv_{+-}(x)-1\big]
\frac{\Ga^{-1}}{\la_-(I_-)\la_+(I_-)}.
\ec
For $x\neq I_-$, these formulas simplify to (\ref{wextr})~(i) and (ii).

Adding formulas (\ref{wextr})~(i) and (iii) yields
\be
w^{(-)}_-(x)+w^{(+)}_-(x)
=\Big[\frac{1}{\la_+(I_-)}+\frac{1}{\la_-(I_+)}\Big]\Vv_{-+}(x)
-\frac{1}{\la_-(I_+)}.
\ee
Since $\Vv_{-+}$ is nonincreasing and continuous, the infimum of this function
over $x\in I$ is equal to the value in $x=I_+$, i.e.,
\be\label{minc}
\inf_{x\in I}\big[w^{(-)}_-(x)+w^{(+)}_-(x)\big]
=\Big[\frac{1}{\la_+(I_-)}+\frac{1}{\la_-(I_+)}\Big]
\frac{\Ga^{-1}}{\la_-(I_+)\la_+(I_+)}
-\frac{1}{\la_-(I_+)}.
\ee
Using the fact that $\Ga>0$, we see that the expression in (\ref{minc}) is
positive if and only if
\be
\frac{1}{\la_+(I_-)\la_+(I_+)}+\frac{1}{\la_-(I_+)\la_+(I_+)}>\Ga.
\ee
Taking into account (\ref{Gamma}) and (\ref{valid})~(ii) (which also holds
with the equality signs reversed), this is equivalent to $f_+(I_+)>0$.
\epro

\subsection{Restricted models}\label{S:restproof}

In the present subsection, we prove Proposition~\ref{P:Phi} as well as
Lemmas~\ref{L:JJ}, \ref{L:comp}, and \ref{L:unif}.\med

\bpro[Proof of Lemma~\ref{L:JJ}]
To prove that $J_-<\phi_+(J_-)$ for all $J_-\in I$, it suffices to show that
$\La_+(J_-,J_-+\eps)>0$ for $\eps>0$ sufficiently small. Here
\be
\La_+(J_-,J_-+\eps)=\frac{1}{\la_+(J_-)\la_+(J_-+\eps)}
+\int_{J_-}^{J_-+\eps}\frac{1}{\la_-}\di\Big(\frac{1}{\la_+}\Big)
\ee
By assumptions (A3) and (A5) and the fact that $J_-\in I$, the first term
tends to a positive limit as $\eps\down 0$ while the second term tends to
zero. By the symmetry between buy and sell orders, we see that also
$\phi_-(J_+)<J_+$ for all $J_+\in I$.

Using (A3), we see that for fixed $J_-$, the function $\La_+(J_-,J_+)$ is
nonincreasing as a function of $J_+$, and hence that $\La_+(J_-,J_+)>0$ if and
only if $J_+<\phi_+(J_-)$. Similarly, $\La_-(J_-,J_+)>0$ if and only if
$\phi_-(J_+)<J_-$, so the second claim of the lemma follows from
Theorem~\ref{T:posrec}, where we use that for $I_-<J_-<J_+<I_+$, the
restricted model on $(J_-,J_+)$ satisfies (A6) since the model on $I$
satisfies (A5).
\epro

\bpro[Proof of Proposition~\ref{P:Phi}]
Recall from Section~\ref{S:restrict} that for any subinterval $J$ such that
$\ov J\sub I$, $\La_\pm(J_-,J_+)$ denote the expressions in (\ref{valid}),
calculated for the process restricted to the subinterval $\ov J$.
We will prove that
\be\label{PsiLa}
V_{\rm W}^{-2}-\Phi(V)
=\La_-\big(\la^{-1}_-(V),\la^{-1}_+(V)\big)
=\La_+\big(\la^{-1}_-(V),\la^{-1}_+(V)\big)
\ee
for any $V\in(V_{\rm W},V_{\rm max}]$ such that
$I_-<\la^{-1}_-(V)<\la^{-1}_+(V)<I_+$.
By Proposition~\ref{P:Luck}, this then implies Proposition~\ref{P:Phi}.

We will first prove (\ref{PsiLa}) for Stigler-Luckock models in standard form
(see Appendix~\ref{A:staform}). Let $D$ be the set of all pairs
$(J_-,J_+)\in\R^2$ such that $I_-\leq J_-\leq J_+\leq I_+$ and
  $\la_-(J_-)=\la_+(J_+)$, and let $T:=\sup\{J_+-J_-:(J_-,J_+)\in D\}$. Note
  that $D\neq\emptyset$ and $T>0$ by (A3) and (A7). We define a curve
$\ga:[0,T]\to D$ with $\ga(t)=\big(\ga_-(t),\ga_+(t)\big)$ by
\be
\ga_-(t):=\inf\big\{J_-\in\ov I:\la_-(t)\leq\la_+(J_-+t)\big\}
\quand\ga_+(t):=\ga_-(t)+t.
\ee
Using (A3), it is not hard to see that $\ga_-$ is nonincreasing, $\ga_+$ is
nondecreasing, and $\ga$ is Lipschitz continuous with Lipschitz constant 1.
Using also (A4), we see that $\ga_-(0)=(x_{\rm W},x_{\rm W})$, where $x_{\rm
  W}$ is the Walrasian price where in fact $x_{\rm W}=0$ since we are
  assuming $\la_\pm$ are in standard form.

For any $J_-,J'_-,J_+,J'_+$ with $J_-\leq J'_-$ and $J_+\leq J'_+$, we observe
that
\bc
\dis\La_+(J'_-,J_+)-\La_+(J_-,J_+)
&=&\dis\int_{J_-}^{J'_-}\Big\{\frac{1}{\la_+(J_+)}-\frac{1}{\la_-}\Big\}
\,\di\Big(\frac{1}{\la_+}\Big),\\[5pt]
\dis\La_+(J_-,J'_+)-\La_+(J_-,J_+)
&=&\dis\int_{J_+}^{J'_+}\Big\{\frac{1}{\la_+(J_-)}+\frac{1}{\la_-}\Big\}
\,\di\Big(\frac{1}{\la_+}\Big).
\ec
Since the model is in standard form, the functions $\la_\pm$ are Lipschitz
continuous with Lipschitz constant one. Since we are also assuming (A5), it
follows that $\la^{-1}_\pm$ are locally Lipschitz on $\{(J_-,J_+):\ov J\sub
I\}$. Therefore, we can apply Lemma~\ref{L:curve} to conclude that
\bc\label{curint}
\dis\La_+\big(\ga(t)\big)
&=&\dis\La_+\big(\ga(0)\big)
+\int_0^t\Big\{\frac{1}{\la_+\circ\ga_+}-\frac{1}{\la_-\circ\ga_-}\Big\}
\,\di\Big(\frac{1}{\la_+\circ\ga_-}\Big)\\[5pt]
&&\dis+\int_0^t\Big\{\frac{1}{\la_+\circ\ga_-}+\frac{1}{\la_-\circ\ga_+}\Big\}
\,\di\Big(\frac{1}{\la_+\circ\ga_+}\Big)
\ec
for any $t\in[0,T')$, where $T':=\sup\{t\in[0,T]:I_-<\ga_-(t)<\ga_+(t)<I_+\}$.

Since $\la_-(\ga_-(s))=\la_+(\ga_+(s))$ for all $s\in[0,t]$, the first
integral in (\ref{curint}) is zero. Set $\phi:=\la_+\circ\ga_+$ and
$\phi^{-1}(V):=\inf\{t\geq 0:\phi(t)\geq V\}$, and observe that
$\ga(\phi^{-1})(V)=\big(\la^{-1}_-(V),\la^{-1}_+(V)\big)$ since the latter is
the smallest interval $J$ such that $\la_-(J_-)=V=\la_+(J_+)$, and $\phi^{-1}$
is left-continuous. Using the substitution of variables $W=\phi(t)$ (recall
(\ref{substitute})), using also the fact that
\be
\La_+\big(\ga(0)\big)=\La_+(x_{\rm W},x_{\rm W})=\frac{1}{V_{\rm W}^2},
\ee
we can rewrite (\ref{curint}) as
\be
\La_+\big(\la^{-1}_-(V),\la^{-1}_+(V)\big)=
\frac{1}{V_{\rm W}^2}
+\int_{V_{\rm W}}^V\Big\{\frac{1}{\la_+\big(\la^{-1}_-(W)\big)}
+\frac{1}{\la_-\big(\la^{-1}_+(W)\big)}\Big\}\di\Big(\frac{1}{W}\Big),
\ee
which holds whevener $V<V_{\rm max}$. This proves (\ref{PsiLa}) for
  $\La_+$. The equality for $\La_-$ follows in the same way, or
alternatively, one can use the fact that
\be\label{Ladif}
\La_+(J_-,J_+)-\La_-(J_-,J_+)
=\Big(\frac{1}{\la_+(J_+)}-\frac{1}{\la_-(J_-)}\Big)
\Big(\frac{1}{\la_+(J_-)}+\frac{1}{\la_-(J_+)}\Big),
\ee
which follows from partial integration of the formulas in (\ref{valid}) and
shows that $\La_+(J_-,J_+)=\La_-(J_-,J_+)$ whenever
$\la_-(J_-)=\la_+(J_+)$.

This proves (\ref{PsiLa}) for models in standard form. Since by
Proposition~\ref{P:stafor} in the appendix, any Stigler-Luckock model can be
brought in standard form, to prove (\ref{PsiLa}) more generally it suffices to
show that the quantities $V_{\rm W}$, $\La_\pm(\la^{-1}_-(V),\la^{-1}_+(V))$,
and $\Phi(V)$ do not change when we bring a model in standard form.

By Proposition~\ref{P:stafor}, for any demand and supply functions $\la'_\pm$
satisfying (A1) and (A2) on some interval $I'$, there exists demand and supply
functions $\la_\pm$ in standard form on some interval $I$, together with a
nondecreasing function $\psi:\ov I\to\ov I'$ that satisfies
$\psi(I_\pm)=I'_\pm$ and that is right-continuous on $I$, such that
$\mu'_\pm=\mu_\pm\circ\psi^{-1}$ and (\ref{ordpres2}) holds.
It is not hard to see that as a result of (\ref{ordpres2}), 
the supply and demand functions $\la_\pm$ and $\la'_\pm$ have the same
Walrasian volume of trade.

In fact, we are only interested here in the case that $\la'_\pm$ satisfy
moreover (A3), (A5), and (A7). In particular, (A3) says that $\mu'_\pm$ have
no atoms in $I'$ and hence $\psi$ is strictly increasing, so (\ref{ordpres2})
is trivial in our setting. Since $\psi$ is strictly increasing, it has an
inverse $\psi^{-1}:\ov I'\to I$ which is a continuous, nondecreasing function.
Since $\mu'_\pm=\mu_\pm\circ\psi^{-1}$, the functions $\la'_\pm$ are given in
terms of $\la_\pm$ by
\be
\la'_\pm(y)=\la_\pm\big(\psi^{-1}(y)\big)\qquad(y\in\ov I').
\ee
Since
\be
\la_\pm\big(\la^{-1}_\mp(V)\big)=\sup\{\la_\pm(x):x\in\ov I,\ \la_\mp(x)\geq V\},
\ee
these quantities are the same for the model in standard form and the
transformed model, and hence the same is true for $\Phi(V)$. By the
substitution of variables formula (see (\ref{substitute})), also
\be
\int_{\la^{-1}_-(V)}^{\la^{-1}_+(V)}\frac{1}{\la_+}\di\Big(\frac{1}{\la_-}\Big)
\quand
\int_{\la^{-1}_-(V)}^{\la^{-1}_+(V)}\frac{1}{\la_-}\di\Big(\frac{1}{\la_+}\Big)
\ee
are the same for the transformed model and hence this is also true for
$\La_\pm(\la^{-1}_-(V),\la^{-1}_+(V))$.
\epro

\bpro[Proof of Lemma~\ref{L:comp}]
Assume that $J\sub I$ satisfies $\ov J\sub I$ and (i) and (ii).
Since $f_-(J_-)=0=f_+(J_+)$, by (\ref{volume}), we
must have $\la_-(J_-)=\la_+(J_+)$. Call this quantity $V$.
If $V=V_{\rm W}$, then it is easy to check that
$\La_\pm(J_-,J_+)=V^{-2}_{\rm W}>0$, which by Proposition~\ref{P:Luck}
contradicts $f_\pm(J_\pm)=0$, so we must have $V>V_{\rm W}$.  Assumption (ii)
implies that $J=\big(\la^{-1}_-(V),\la^{-1}_+(V)\big)$. Now
Proposition~\ref{P:Phi} tells us that $\Phi(V)=V_{\rm W}^{-2}$ and hence
$V=V_{\rm L}$ and $J=J^\cri$.

Conversely, if $V_{\rm L}<V_{\rm max}$, then $J^\cri:=\big(\la^{-1}_-(V_{\rm
  L}),\la^{-1}_+(V_{\rm L})\big)$ satisfies $\ov J\sub I$ and
Proposition~\ref{P:Phi} tells us that the solution
to Luckock's equation on $\ov J$ satisfies $f_-(J_-)=0=f_+(J_+)$
while (ii) holds because of the way $\la^{-1}_\pm$ have been defined.
\epro

\bpro[Proof of Lemma~\ref{L:unif}]
For the uniform model with $\la_-(x)=1-x$ and $\la_+(x)=x$, we have
\be
V_{\rm W}=\ha,\quad V_{\rm max}=1,\quad
\la^{-1}_-(V)=1-V,\quand\la^{-1}_+(V)=V.
\ee
It follows that the function $\Phi$ from (\ref{Phidef}) is given by
\be\ba{l}
\dis\Phi(V)=-2\int_{1/2}^V\frac{1}{1-W}\di\Big(\frac{1}{W}\Big)
=2\int_{1/V}^2\frac{1}{1-\frac{1}{y}}\,\di y\\[5pt]
\dis\quad=2\int_{1/V}^2\big\{1+\frac{1}{y-1}\big\}\,\di y
=4-2\big\{V^{-1}+\log(V^{-1}-1)\big\}.
\ec
Setting $\Phi(V)=V_{\rm W}^{-2}$ gives
\be
-V^{-1}=\log\big(V^{-1}-1\big)
\quad\desd\quad
\ex{-V^{-1}}=V^{-1}-1,
\ee
which tells us that $V_{\rm L}=1/z$ where $z$ solves
$f(z):=e^{-z}-z+1=0$. Since the function $f$ is continuous and strictly
decreasing with $f(1)=e^{-1}$ and $f(z)\to-\infty$ for $z\to\infty$, the
equation $f(z)=0$ has a unique solution $z$, and this solution satisfies
$z>1$.
\epro


\section{Analysis of the Markov chain}\label{S:Markov}

\subsection{A consequence of stationarity}\label{S:station}

In this subsection, we prove Theorem~\ref{T:Luckock}.\med

\bpro[Proof of Theorem~\ref{T:Luckock}]
Let $\nu$ be an invariant law, let $X_0$ be an $\Si_{\rm ord}$-valued
random variable with law $\nu$, and let $(U_1,\sig_1)$ be independent of $X_0$
with law $\ov\mu:=|\mu|^{-1}\mu$, as defined in Section~\ref{S:intro}. Then
stationarity means that $X_1:=L_{U_1,\sig_1}(X_0)$ has the same law as $X_0$,
where $L_{u,\sig}$ is the Luckock map from (\ref{Luckmap}).

Set $M_\pm:=M_\pm(X_0)$ and let $\mu_\pm$ be as in (\ref{mupm}). We claim that
\be\ba{rr@{\,}c@{\,}l}\label{inteq}
{\rm(i)}&\dis\int_A\P[M_-<x]\,\mu_+(\di x)
&=&\dis\int_A\la_-(x)\,\P[M_+\in\di x],\\[7pt]
{\rm(ii)}&\dis\int_A\P[M_+>x]\,\mu_-(\di x)
&=&\dis\int_A\la_+(x)\,\P[M_-\in\di x]
\ec
for each measurable $A\sub I$ that is contained in some compact subinterval
$[J_-,J_+]\sub I$. Indeed, stationarity implies (see \cite[Lemma~10]{Swa15})
that for any measurable $A\sub I$, sell limit orders are added in $A$ with the
same frequency as they are removed, i.e.,
\be\label{addremove}
\P\big[X_1^+(A)=X_0^+(A)+1\big]
=\P\big[X_1^+(A)=X_0^+(A)-1\big].
\ee
Recalling the definition of the Luckock map in (\ref{Luckmap}), we see that
this means that
\be
\P\big[\sig_1=+,\ U_1\in A,\ M_-<U_1\big]
=\P\big[\sig_1=-,\ M_+\in A,\ M_+\leq U_1\big].
\ee
Since $(U_1,\sig_1)$ has law $\ov\mu$ and is independent of $M_\pm$, it
follows that
\be
|\mu|^{-1}\int_A\mu_+(\di x)\P[M_-<x]
=|\mu|^{-1}\int_A\P[M_+\in\di x]\mu_-\big([x,I_+]\big),
\ee
which up to the factor $|\mu|^{-1}$ is (\ref{inteq})~(i).
Similarly, equation (\ref{inteq})~(ii) follows from the requirement that 
buy limit orders are added in $A$ with the same frequency as they
are removed.

By (A3), the measures $\mu_\pm$ do not have atoms in $I$ and hence by
(A5) and (\ref{inteq}), the same is true for the laws of $M_\pm$. It follows
that $\P[M_-<x]=f_-(x)$ and $\P[M_+>x]=f_+(x)$ $(x\in I)$ and $f_\pm$ are
continuous.  Since by (\ref{mupm}), $-\di\la_-$ and $\di\la_+$ are the
restrictions of the measures $\mu_-$ and $\mu_+$ to $I$, respectively, we can
rewrite (\ref{inteq})~(i) and (ii) as (\ref{Luckock})~(i) and (ii). Since
$f_-(I_+)=\P[M_-\leq I_+]=1$ and $f_+(I_-)=\P[M_+\geq I_-]=1$, the boundary
conditions (\ref{Luckock})~(iii) also follow.
\epro

\subsection{A Lyapunov function}

It follows from Theorem~\ref{T:Luckock} that if a Stigler-Luckock model is
positive recurrent, then the solution to Luckock's equation must satisfy
$f_-(I_-)\wedge f_+(I_+)>0$. Theorem~\ref{T:posrec} states that this condition
is also sufficient. We will prove this by showing that
\be\label{Ldef}
\Ly(\Xc):=\sqrt{(F^{(-)}(\Xc)\vee 0)^2+(F^{(+)}(\Xc)\vee 0)^2}
\qquad\big(\Xc\in\Si^{\rm fin}_{\rm ord}\big)
\ee
is a Lyapunov function. We note that this is the only place in the paper where
we make use of a function of a Stigler-Luckock process that is not linear
(namely $\Ly$). In view of Theorem~\ref{T:inverse}, we have fairly good
control of linear functionals, which as in Theorem~\ref{T:Luckock} (which
depends on Theorem~\ref{T:spec}) allows us to more or less explicitly
calculate the marginal distributions of the best buy and sell offers
$M_-(\Xc)$ and $M_+(\Xc)$ in equilibrium. Proving that a Stigler-Luckock model
is positive recurrent, however, always entails proving something about the
joint distribution of $M_-(\Xc)$ and $M_+(\Xc)$. Indeed, the following
proposition can be used to give a lower bound on the probability, in
equilibrium, that the order book is empty, which corresponds to the event that
$M_-(\Xc)=I_-$ and at the same time $M_+(\Xc)=I_+$, but this bound is not very
explicit or sharp. It seems that such information cannot be obtained
from linear functionals and indeed for no choice of weight functions
$(w_-,w_+)$ is a linear function of the form (\ref{Flin}) a Lyapunov function.

Recall the definition (\ref{Gdef}) of the generator $G$ of a Stigler-Luckock
model. We have the following result.

\bp[Lyapunov function]
Assume\label{P:Lyap} (A3) and (A6), and that the unique solution $(f_-,f_+)$ of Luckock's
equation (\ref{Luckock}) satisfies $\eps:=f_-(I_-)\wedge f_+(I_+)>0$. Then
there exists a constant $K<\infty$ such that the function in (\ref{Ldef})
satisfies $G\Ly(\Xc)\leq K$ for all $\Xc\in\Si^{\rm fin}_{\rm ord}$. Moreover,
for each $\eps'<\eps$, there exists an $N<\infty$ such that
\be\label{GLya}
G\Ly(\Xc)\leq-\eps'\quad\mbox{whenever}\quad|\Xc^-|+|\Xc^+|\geq N.
\ee
\ep
\bpro
Let us write $\vec F(\Xc):=\big(F^{(-)}(\Xc),F^{(+)}(\Xc)\big)$ and let
$|\,\cdot\,|$ denote the euclidean norm on $\R^2$. Let us also
write $\Ly(\Xc)=v(\vec F(\Xc))$ where $v:R^2\to\R$ is the function
\be
v(z_1,z_2):=\sqrt{(z_1\vee 0)^2+(z_2\vee 0)^2}.
\ee
Set
\be\label{Wdef}
W:=\sup_{x\in I}\big|\big(w^{(-)}_-(x),w^{(+)}_-(x)\big)\big|
\vee\sup_{x\in I}\big|\big(w^{(-)}_+(x),w^{(+)}_+(x)\big)\big|,
\ee
which is the maximal amount by which $\vec F(\Xc)$ can change due to the
addition or removal of a single limit order. Since the function $v$ is
Lipschitz continuous with Lipschitz constant 1, we can estimate
\bc
\dis G\Ly(\Xc)
&=&\dis\int\big\{\Ly\big(L_{u,\sig}(\Xc)\big)-\Ly\big(\Xc\big)\big\}
\,\mu\big(\di(u,\sig)\big)\\[5pt]
&\leq&\dis\int\big|\vec F\big(L_{u,\sig}(\Xc)\big)-\vec F\big(\Xc\big)\big|
\,\mu\big(\di(u,\sig)\big)\\[5pt]
&\leq&\dis W\big(\la_-(I_-)+\la_+(I_+)\big)=:K.
\ec
Let
\be
\de:=\inf_{x\in I}\big[w^{(-)}_-(x)+w^{(+)}_-(x)\big]
\wedge\inf_{x\in I}\big[w^{(-)}_+(x)+w^{(+)}_+(x)\big],
\ee
which is positive by (\ref{wsum}) and our assumption that
$f_-(I_-)\wedge f_+(I_+)>0$. Since adding a limit order to the order book
always raises $F^{(-)}+F^{(+)}$ by at least $\de$,
\be\label{delow}
F^{(-)}(\Xc)+F^{(+)}(\Xc)\geq\de\big(|\Xc^-|+|\Xc^+|\big).
\ee
This shows that $\vec F(\Xc)$ takes values in the half space
$H:=\{(z_1,z_2)\in\R^2:z_1+z_2>0\}$ as long as $\Xc\neq 0$, and moreover
$\big|\vec F(\Xc)\big|$ is large if $|\Xc^-|+|\Xc^+|$ is.

For any $z=(z_1,z_2)\in\R^2$ with $z_1+z_2>0$, let us define
\be
p_1(z):=\frac{z_1\vee 0}{\sqrt{(z_1\vee 0)^2+(z_2\vee 0)^2}}
\quand
p_2(z):=\frac{z_2\vee 0}{\sqrt{(z_1\vee 0)^2+(z_2\vee 0)^2}}.
\ee
Then, for any $y,z\in H$, we can write
\be
v(z)=v(y)+p_1(y)(z_1-y_1)+p_2(y)(z_2-y_2)+R(y,z),
\ee
where for any $y,z$ that differ at most by the constant $W$ from (\ref{Wdef}),
the error term $R(x,y)$ can be estimated as
\be
R(y,z)\leq C|y|^{-1}\qquad\big(y,z\in H,\ |z-y|\leq W\big)
\ee
for some constant $C<\infty$. It follows that we can write
\be\label{GLy}
G\Ly(\Xc)=p_1(\vec F(\Xc))GF^{(-)}(\Xc)+p_2(\vec F(\Xc))GF^{(+)}(\Xc)
+E(\Xc),
\ee
where the error term can be estimated as
\bc\label{Ebd}
\dis\big|E(\Xc)|
&=&\dis\Big|\int R\big(\vec F(\Xc),\vec F(L_{u,\sig}(\Xc)\big)
\,\mu\big(\di(u,\sig)\big)\Big|\\[5pt]
&\leq&\dis C\big(\la_-(I_-)+\la_+(I_+)\big)\big|\vec F(\Xc)\big|^{-1},
\ec
which in view of (\ref{delow}) can be made arbitrary small by choosing 
$|\Xc^-|+|\Xc^+|$ sufficiently large.

By Theorem~\ref{T:spec} and the way we have defined the weight functions
$w^{(-)}_\pm$ and $w^{(+)}_\pm$, one has
\bc\label{Gifpos}
\dis GF^{(-)}(\Xc)&=&\dis-f_-(I_-)\mbox{ if }|\Xc^-|\neq 0,\\[5pt]
\dis GF^{(+)}(\Xc)&=&\dis-f_+(I_+)\mbox{ if }|\Xc^+|\neq 0.
\ec
It follows from (\ref{wextr}) and elementary properties of the functions in
(\ref{Vpm}) that
\be
w^{(-)}_->0,\quad w^{(-)}_+<0,\quad w^{(+)}_-<0,\quand
w^{(+)}_+>0\quad\mbox{on }I.
\ee
In view of this, we have
\be\ba{r@{\;\;}c@{\;\;}l@{\;\;}c@{\;\;}l}
\dis|\Xc^-|=0&\dis\volgt&\dis F^{(-)}(\Xc)\leq 0
&\dis\volgt&\dis p_1(\vec F(\Xc))=0,\\[5pt]
\dis|\Xc^+|=0&\dis\volgt&\dis F^{(+)}(\Xc)\leq 0
&\dis\volgt&\dis p_2(\vec F(\Xc))=0.
\ec
Combining this with (\ref{Gifpos}), we obtain that
\be
p_1(\vec F(\Xc))GF^{(-)}(\Xc)+p_2(\vec F(\Xc))GF^{(+)}(\Xc)\leq-\eps
\qquad(\Xc\neq 0).
\ee
Inserting this into (\ref{GLy}), using our bound (\ref{Ebd}) on the error
term, and using also (\ref{delow}), we see that by choosing $|\Xc^-|+|\Xc^+|$
large enough, we can make $G\Ly$ smaller than $-\eps'$ for any $\eps'<\eps$. 
\epro

\subsection{Positive recurrence}\label{S:posprf}

\bpro[Proof of Theorem~\ref{T:posrec}]
If a Stigler-Luckock model is positive recurrent, then it it is possible to
construct a stationary process $(X_k)_{k\in\Z}$ that makes i.i.d.\ excursions
away from the empty state $0$. In particular, positive recurrence implies the
existence of an invariant law $\nu$ that is concentrated on $\Si^{\rm
  fin}_{\rm ord}$ and satisfies $\nu(\{0\})>0$. By Theorem~\ref{T:Luckock}, it
follows that Luckock's equation (\ref{Luckock}) has a solution $(f_-,f_+)$
such that $f_-(I_-)\wedge f_+(I_+)\geq\nu(\{0\})>0$.

Conversely, assume (A3) and (A6) and that the (by Proposition~\ref{P:Luck}
unique) solution to Luckock's equation satisfies $f_-(I_-)\wedge f_+(I_+)>0$.
Let $P$ denote the transition kernel of the discrete-time process
$(X_k)_{k\geq 0}$ and for any nonnegative measurable function $f:\Si^{\rm
  fin}_{\rm ord}\to\R$ write $Pf(x):=\int P(x,\di y)f(y)$. Write
\be
C_N:=\{\Xc\in\Si^{\rm fin}_{\rm ord}:|\Xc|^-+|\Xc|^+<N\}.
\ee
Multiplying the Lyapunov function $V$ of Proposition~\ref{P:Lyap} by a
suitable constant, we obtain a nonnegative function $f$ and finite constants
$K,N$ such that
\be
Pf-f\leq K1_{C_N}-1.
\ee
Let $\tau_0:=\inf\{k>0:X_k=0\}$ denote the first return time to the empty
configuration. By assumption (A6), there exists a constant $\eps>0$ such that
\be
\P^x[\tau_0\leq N+1]\geq\eps\qquad(x\in C_N).
\ee
Moreover, (A6) guarantees that $\P^0[\tau_0=0]>0$, which shows that the model is
aperiodic from $0$. Applying Proposition~\ref{P:posrec} from
Appendix~\ref{A:posrec} in the appendix, we conclude that the Stigler-Luckock
model under consideration is positive recurrent and (\ref{tonu}) holds.
\epro

\appendix

\section{Appendix}

\subsection{The model in standard form}\label{A:staform}

By definition, we say that a Stigler-Luckock model is in \emph{standard form}
if
\be\label{standard}
\la_+(x)-\la_-(x)=x\qquad(x\in\ov I).
\ee
Note that for a model in standard form, $\di\la_+-\di\la_-$, which is the
total rate at which limit orders arrive in $I$, is the Lebesgue measure.
Clearly, a model in standard form always satisfies (A3) and (A4). Let
$p_+:\ov I\to[0,1]$ be the Radon-Nikodym derivative of $\di\la_+$ w.r.t.\ the
Lebesgue measure and let $p_-:=1-p_+$. Then
\be
\la_-(x)=\la_-(I_+)+\int_x^{I_+}\! p_-(x)\,\di x
\quand
\la_+(x)=\la_+(I_-)+\int_{I_-}^x\! p_+(x)\,\di x.
\ee
In particular, a model in standard form is uniquely characterized by its
interval $I$, rates of market makers $\la_\pm(I_\mp)$, and the function $p_+$.

In the present section we will show that under mild additional assumptions,
each Stigler-Luckock model satisfying (A1) and (A2) can be brought in standard
form. To demonstrate the idea on a concrete example, consider a
Stigler-Luckock model $(\Xc_t)_{t\geq 0}$ with $I=(0,2n)$ where $n\geq 1$ is
some integer, there are no market orders, and
\be
p_+(x)=1_{\{\lfloor x\rfloor\mbox{ is even}\}}\qquad(x\in I).
\ee
Assume that $\Xc_0=0$ (the empty initial state). Define a function $\psi: I\to
I':=(0,n+1)$ by $\psi(x):=1+\lfloor x/2\rfloor$, and let
$\Xc'_t:=\Xc_t\circ\psi^{-1}$ $(t\geq 0)$ denote the image of $\Xc_t$ under
$\psi$. Then $(\Xc'_t)_{t\geq 0}$ is a Stigler-Luckock model where limit
orders are placed at discrete prices only. More precisely, for each
$k\in\{1,\ldots,n\}$, in the original model, buy and sell limit orders arrive
in $(2k-2,2k)$ with rate one each, in such a way that sell limit orders arrive
on the left of buy limit orders. After applying the map $\psi$, all these
orders arrive at the price $k$, where they still match.

Our aim is to show that this construction works quite generally, i.e., given a
fairly general Stigler-Luckock model $(\Xc'_t)_{t\geq 0}$ with parameters $I'$
and $\la'_\pm$, we can find a Stigler-Luckock model $(\Xc_t)_{t\geq 0}$ in
standard form with parameters $I$ and $\la_\pm$, as well as a nondecreasing
right-continuous function $\psi:I\to I'$, such that
$\Xc'_t=\Xc_t\circ\psi^{-1}$ $(t\geq 0)$. We will see that moreover, the
parameters of $(\Xc'_t)_{t\geq 0}$ determine those of the model in standard
form as well as the function $\psi$ uniquely.

We first investigate which functions $\psi$ have the property that the image
of a Stigler-Luckock model under $\psi$ is again a Stigler-Luckock
model. Clearly, it suffices if $\psi$ is strictly increasing, but as our
previous example showed, this condition can be weakened.

\bl[Transformed Stigler-Luckock model]\hspace*{1pt}
Let\label{L:ordpres} $(\Xc_t)_{t\geq 0}$ be a Stigler-Luckock model on some
nonempty open interval $I\sub\R$ with demand and supply functions $\la_\pm$
satisfying (A1) and (A2). Let $I'\sub\R$ be an open interval, let $\psi:I\to
I'$ be nonincreasing, and assume that
\be\label{ordpres}
\int_I\!\nu_-(\di x)\int_I\!\nu_+(\di y)\,
1_{\{x<y\mbox{ and }\psi(x)\geq\psi(y)\}}=0,
\ee
where $\nu_\pm:=\mu_\pm+\Xc^\pm_0$ with $\mu_\pm$ defined in (\ref{mupm}).
Extend $\psi$ to $\ov I$ by setting $\psi(I_\pm):=I'_\pm$. Then
$\Xc'_t:=\Xc_t\circ\psi^{-1}$ $(t\geq 0)$ is a Stigler-Luckock model with
demand and supply functions $\la'_\pm$ and related measures $\mu'_\pm$ as in
(\ref{mupm}) given by $\mu'_\pm:=\mu_\pm\circ\psi^{-1}$.
\el
\bpro
The condition (\ref{ordpres}) guarantees that a.s., of all the orders that are
at time zero in the order book or that arrive at later times, a buy and a sell
order match in the original model $\Xc$ if and only if they match in the
transformed model $\Xc'$.
\epro

We wish to show that a general Stigler-Luckock model $(\Xc'_t)_{t\geq 0}$ can
be obtained as a function of a Stigler-Luckock model $(\Xc_t)_{t\geq 0}$ in
standard form. We will need a weak assumption on the initial state of
$(\Xc'_t)_{t\geq 0}$. To formulate this properly, we need some definitions.
Let $\mu$ be a finite nonnegative measure on $\ov\R:=[-\infty,\infty]$ and let
${\rm supp}(\mu)$ denote its support. Then the complement of ${\rm supp}(\mu)$
is a countable union of disjoint open intervals. If for each left endpoint
$x_-$ of such an interval $(x_-,x_+)$, we remove $x_-$ from ${\rm
  supp}(\mu)$ if it carries no mass, then we obtain
\be
{\rm supp}_+(\mu)=
\big\{x\in\ov\R:\mu\big([x,y)\big)>0\ \forall y>x\big\}.
\ee
The set ${\rm supp}_+(\mu)$ is the support of $\mu$ with respect to the
topology of convergence from the right, where a sequence $x_n$ converges to a
limit $x$ if and only if $x_n\to x$ in the usual topology on $\ov\R$ and
moreover $x_n\geq x$ for $n$ large enough. A basis for this topology is formed
by all sets of the form $[x,y)$ with $x<y$.

\bp[Standard form]
Let\label{P:stafor} $(\Xc'_t)_{t\geq 0}$ be a Stigler-Luckock model with
demand and supply functions $\la'_\pm:\ov I'\to\half$ satisfying (A1) and
(A2). Assume that $\la'_\pm$ are not both constant and that
\be\label{initcond}
\Xc'_0\mbox{ is concentrated on }{\rm supp}_+(\mu'_-+\mu'_+),
\ee
where $\mu'_\pm$ are defined in terms of $\la'_\pm$ as in (\ref{mupm}). Then
there exist demand and supply functions $\la_\pm$ in standard form on some
interval $I$, as well as a nondecreasing function $\psi:\ov I\to\ov I'$ that
satisfies $\psi(I_\pm)=I'_\pm$ and that is right-continuous on $I$, such that 
$\mu'_\pm=\mu_\pm\circ\psi^{-1}$ and
\be\label{ordpres2}
\int_I\!\mu_-(\di x)\int_I\!\mu_+(\di y)\,
1_{\{x<y\mbox{ and }\psi(x)\geq\psi(y)\}}=0.
\ee
Moreover, these conditions determine $I_\pm$, $\la_\pm$, and $\psi$
uniquely. Finally, there exists a Stigler-Luckock model $(\Xc_t)_{t\geq 0}$
with demand and supply functions $\la_\pm$ such that (\ref{ordpres}) holds and
the process $\Xc''_t:=\Xc_t\circ\psi^{-1}$ $(t\geq 0)$ is equal in law to
$(\Xc'_t)_{t\geq 0}$.
\ep
\bpro
Let $\chi(y):=\la'_+(y)-\la'_-(y)$ $(y\in\ov I)$. Since $\la'_\pm$ satisfy (A1)
and (A2), the function $\chi$ is nondecreasing on $\ov I'$, and continuous at
$I'_\pm$. We claim that our conditions imply that we must choose
\be\label{If}
{\rm(i)}\quad I_\pm=\chi(I'_\pm)
\qquad{\rm(ii)}\quad\psi(x)=\sup\{y\in I':\chi(y)\leq x\}\quad(x\in I).
\ee
Indeed, $\chi(I'_-)=\la'_+(I'_-)-\la'_-(I'_-)$ is the rate of sell market orders
minus the total rate of buy orders for the process $\Xc'$. Since
$\mu'_\pm=\mu_\pm\circ\psi^{-1}$, this must equal the same quantity for the
process $\Xc$, i.e., $\la_+(I_-)-\la_-(I_-)$, which equals $I_-$ by the
assumption that $\Xc$ is in standard form. The same argument shows that
$I_+=\chi(I'_+)$. Let $\rho$ be the restriction to $I$ of $\mu_++\mu_-$ and let
$\rho'$ be similarly defined for the process $\Xc'$. Then $\rho$ is the
Lebesgue measure on $I$ by the assumption that $\Xc$ is in standard form.
Define $\psi^{-1}:I'\to[I_-,I_+)$ by
\be
\psi^{-1}(y):=\inf\{x\in I:\psi(x)\geq y\}\qquad(y\in I'),
\ee
where the infimum of the empty set is $:=I_-$. Then, for $y\in I'$,
\be\label{invint}
\psi^{-1}\big((I'_-,y)\big)
:=\big\{x\in I:\psi(x)\in(I'_-,y)\big\}
=\big(I_-,\psi^{-1}(y)\big).
\ee
Letting $\chi(y-)$ denote the left-continuous modification of $\chi$,
using (\ref{If})~(i) and the fact that $\rho'_\pm=\rho_\pm\circ\psi^{-1}$,
it follows that
\be\label{chinv}
\chi(y-)=\chi(I'_-)+\rho'\big((I'_-,y)\big)
=I_-+\rho\big((I'_-,\psi^{-1}(y))\big)
=\psi^{-1}(y),
\ee
where in the last step we have used that $\rho$ is the Legesgue measure.
This proves that $\chi(y-)$ is the left-continuous inverse of $\psi$ and hence
that $\psi$ is the right-continuous inverse of $\chi$, completing the proof of
(\ref{If})~(ii).

In particular, this shows that $I$ and $\psi$ are uniquely determined by our
conditions. Conversely, choosing $I$ and $\psi$ as in (\ref{If}), our
arguments show that $\rho'=\rho\circ\psi^{-1}$ and that $\Xc$ and $\Xc'$ have
the same rates of buy and sell market orders. Let $\rho_\pm$
denote the restrictions of $\mu_\pm$ to $I$, and let $\rho'_\pm$ be defined
similarly for the process $\Xc'$. To complete the proof, it suffices to show
that we can satisfy $\rho'_\pm=\rho_\pm\circ\psi^{-1}$ and (\ref{ordpres2}).
Let $p_\pm:I\to[0,1]$ and $p'_\pm:I'\to[0,1]$ be defined as the Radon-Nikodym
derivatives
\be
p_\pm:=\frac{\di\mu_\pm}{\di\rho}
\quand
p'_\pm:=\frac{\di\mu'_\pm}{\di\rho'}.
\ee
To specify $\rho_\pm$ in terms of $\rho$, it suffices to specify $p_\pm$.
Letting $\chi(y-)$ and $\chi(y+)$ denote the left- and right-continuous
modifications of $\chi$, we have, for $y\in I$,
\be\label{fjumps}
\chi(y+)-\chi(y)=\mu'_-\big(\{y\})
\quand
\chi(y)-\chi(y-)=\mu'_+\big(\{y\}).
\ee
We will choose
\be
p_\pm(x):=p'_\pm(y)\qquad\mbox{if }\psi(x)=y,\ \rho(\{y\})=0,
\ee
and
\be\left.\bac\label{pmatch}
\dis p_-(x)&:=&\dis1_{(\chi(y),\chi(y+))}\\[5pt]
\dis p_+(x)&:=&\dis1_{(\chi(y-),\chi(y))}
\ea\right\}\qquad\mbox{if }\psi(x)=y,\ \rho(\{y\})>0.
\ee
In view of (\ref{fjumps}), choosing $p_\pm$ in this way ensures that
$\rho'_\pm=\rho_\pm\circ\psi^{-1}$. Moreover, (\ref{pmatch}) guarantees that
if $\psi$ is constant on an interval of the form $[\chi(y-),\chi(y_+))$, then
buy orders are placed on the right of sell orders, so that (\ref{ordpres2}) is
satisfied. It is not hard to see that this is the only way to choose
$\rho_\pm$ such that $\rho'_\pm=\rho_\pm\circ\psi^{-1}$ and (\ref{ordpres2})
is satisfied.

Note that (\ref{ordpres}) simplifies to (\ref{ordpres2}) if $\Xc_0$ is the
empty initial state. In view of this and Lemma~\ref{L:ordpres}, our proof is
complete if $\Xc'$ is started in the empty initial state. To prove the
statement for more general initial states, we must show that if
(\ref{initcond}) holds, then we can choose $\Xc_0$ such that
$\Xc_0\circ\psi^{-1}=\Xc'_0$ and (\ref{ordpres}) is satisfied. We choose
\be
\Xc^\pm_0:=\Xc^{\pm\,'}_0\circ\chi_\mp^{-1},
\ee
where $\chi_-$ and $\chi_+$ denote the left- and right-continuous
modifications of the function $\chi$, respectively, i.e., a buy (resp.\ sell)
limit order of $\Xc'_0$ at the price $y$ is represented in $\Xc_0$ by a buy
(resp.\ sell) limit order at $\chi_+(y)$ (resp.\ $\chi_-(y)$).

Let $\psi(I)$ denote the image of $I$ under $\psi$.
Since $\chi_\pm$ are the left- and right-continuous inverses of $\psi$, we
have $\chi_\pm(y)\in I$ and $\psi(\chi_\pm(y))=y$ for all $y\in\psi(I)$
We claim that $\psi(I)={\rm supp}_+(\mu'_-+\mu'_+)\cap I'$. Indeed, for any
$y_1,y_2\in I'$ with $y_1<y_2$, one has $\psi(I)\cap[y_1,y_2)=\emptyset$ if
and only if $\chi_-(y_1)=\chi_-(y_2)$. By (\ref{invint}) and (\ref{chinv}),
\be
\rho'\big([y_1,y_2)\big)=\rho\big([\chi(y_1-),\chi(y_2-))\big)
=\chi(y_2-)-\chi(y_1-),
\ee
so $\psi(I)\cap[y_1,y_2)=\emptyset$ if and only if $\rho\big([y_1,y_2)\big)=0$.
In view of this, (\ref{initcond}) guarantees that $\Xc'_0=\Xc_0\circ\psi^{-1}$.
(Indeed, the only reason why we need (\ref{initcond}) is that the image of
$\Xc_0$ under $\psi$ can only contain limit orders in $\psi(I)$, which is
${\rm supp}_+(\mu'_-+\mu'_+)\cap I'$.)

To complete the proof, we must show that (\ref{ordpres2}) and our choice of
$\Xc_0$ imply (\ref{ordpres}). The only way (\ref{ordpres}) can fail is that
there exists some $y\in I'$ such that $\chi_-(y)<\chi_+(y)$ and hence $\psi$
is constant on the interval $J:=[\chi_-(y),\chi_+(y))$, while
\be\label{ordpres3}
\int_J\!\nu_-(\di x)
\int_J\!\nu_+(\di y)1_{\{x<y\}}>0.
\ee
Since we represent a buy (resp.\ sell) limit order of $\Xc'_0$ at the price
$y$ by a buy (resp.\ sell) limit order of $\Xc_0$ at $\chi_+(y)$
(resp.\ $\chi_-(y)$), the only atoms of $\Xc^+_0$ in the interval
$J$ must be located at $\chi_-(y)$ while $\Xc^-_0$ cannot
have atoms in $J$. In view of this, (\ref{ordpres2}) implies that
(\ref{ordpres3}) cannot hold.
\epro

\subsection{Ergodicity of Markov chains}\label{A:posrec}

Let $(E,\Ei)$ be a measurable space and let $P$ be a measurable probability
kernel on $E$. For simplicity, we assume that the one-point sets are
measurable, i.e., $\{x\}\in\Ei$ for all $x\in E$. It is known
\cite[Thm~3.4.1]{MT09} that for each probability measure $\mu$ on $E$ there
exists a Markov chain $X=(X_k)_{k\geq 0}$, unique in distribution, such that
$X_0$ has law $\mu$ and the conditional law of $X_{k+1}$ given
$(X_0,\ldots,X_k)$ is given by $P$, for each $k\geq 0$. (This statement is
not quite as straightforward as it may sound since for general measurable
spaces, Kolmogorov's extension theorem is not available.)

We let $P^k$ denote the $k$-th power of $P$. For any measurable real function
$f:E\to[-\infty,\infty]$, we write $P^kf(x):=\int_EP^k(x,\di y)f(y)$, as long
as the integral is well-defined for all $x\in E$. For any probability measure
$\mu$ on $E$ we let $\mu P^k(A):=\int_E\mu(\di x)P^k(x,A)$ $(A\in\Ei)$. Then
$\mu P^k$ is the law of $X_k$ if $X_0$ has law $\mu$. An invariant law of $X$
is a probability measure $\nu$ such that $\nu P=\nu$. We let $\|\mu-\nu\|$
denote the total variation norm distance between two probability measures
$\mu$ and $\nu$.

For any point $x\in E$, let $\P^x$ denote the law of the Markov chain $X$
started from $X_0=x$. Let $\tau_x:=\inf\{k>0:X_k=x\}$ denote the first return
time to $x$. We say that the Markov chain $X$ is \emph{aperiodic from $x$}
if the greatest common divisor of $\{k>0:\P^x[\tau_x=k]>0\}$ is one.

Markov chains satisfying the conditions (\ref{Lyap}) and (\ref{atom}) below
behave in many ways like positive recurrent Markov chains with countable state
space. In particular, (\ref{Lyap}) says that $f$ is a Lyapunov function that
guarantees that the return times to the set $C$ have finite expectation, while
(\ref{atom}) says that once the chain enters $C$, there is a uniformly
positive probability of entering the atom $0$ after a finite number of steps.

\bp[Ergodicity for positive point recurrent chain]
Fix\label{P:posrec} a point $0\in E$. Assume that there exists a measurable function
$f:E\to[0,\infty)$, a measurable set $0\in C\sub E$, and constants $F,K<\infty$
such that $\sup_{x\in C}f(x)\leq F$ and
\be\label{Lyap}
Pf-f\leq K1_C-1.
\ee
Assume moreover that there exist constants $\eps>0$ and $k\geq 0$ such that
\be\label{atom}
\P^x[\tau_0\leq k]\geq\eps\qquad(x\in C).
\ee
Then $\E^x[\tau_0]<\infty$ for all $x\in E$, and the Markov chain $X$ has a
unique invariant law $\nu$. If moreover $X$ is aperiodic from $0$, then
\be\label{ergo}
\|\mu P^n-\nu\|\asto{n}0
\ee
for each probability measure $\mu$ on $E$.
\ep
\bpro
Let $\tau_C:=\inf\{k\geq 1:X_k\in C\}$ denote the first entry time of $C$.
Then \cite[Thm~11.3.4]{MT09} tells us that $\E^x[\tau_C]\leq f(x)+K1_C(x)$
($x\in E$). Since after each visit to $C$, by (\ref{atom}) there is a
probability of at least $\eps$ to visit $0$ in the next $k$ steps, it is not
hard to deduce that $\E^x[\tau_0]<\infty$ for all $x\in E$. Again by
\cite[Thm~11.3.4]{MT09} and the fact that, in the light of (\ref{atom}), $C$
is petite as defined in \cite[Section~5.5.2]{MT09}, we have that $X$ is
positive Harris recurrent. In particular, by \cite[Thm~10.0.1]{MT09}, $X$ has a
unique invariant law $\nu$. Since $\E^x[\tau_0]<\infty$ for all $x\in E$, it
is easy to see that $\nu(\{0\})>0$. By \cite[Thm~10.4.9]{MT09}, $\nu$ is
equivalent to the measure $\psi$ from \cite[Prop.~4.2.2]{MT09}, so
aperiodicity from $0$ as we have defined it implies $\psi$-aperiodicity as
defined in \cite[Section~5.4.3]{MT09}. Now (\ref{ergo}) follows from
\cite[Thm~13.3.3]{MT09}.
\epro

\subsection{Discrete models}\label{A:dis}

Often, it is natural to consider Stigler-Luckock models where the interval $I$
is of the form $I=[0,n]$, with $n\geq 2$ an integer, and the measures
$\mu_\pm$ that determine the rate at which orders arrive are supported on the
set of integers $\{0,\ldots,n\}$. One motivation for this is that real prices
take values that differ by a minimal amount, the so called tick size. Also,
the numerical data for the uniform model shown in Figures~\ref{fig:R} and
\ref{fig:numer} are based on approximation with discrete models with a high
value of $n$.

Although, in the light of Appendix~\ref{A:staform}, discrete Stigler-Luckock
models are in principle included in our general analysis, in practice, when
doing (numerical) calculations, it is more convenient to replace the
differential equations for the general model by difference equations. It turns
out that these difference equations can be solved explicitly much in the
same way as the differential equations of the general model.

In the discrete setting, it is convenient to reparametrize the model
somewhat. We replace the set $\{1,2,\ldots,n\}$ of possible prices of buy
orders by $\{4,6,\ldots,2n+2\}$ and we let $2$ (instead of $0$) be the value
of $M_-(\Xc)$ that signifies that the order book contains no buy limit orders.
Likewise, for sell orders or $M_+(\Xc)$, we replace the set of possible prices
$\{0,1,\ldots,n\}$ by $\{1,3,\ldots,2n+1\}$. Note that in this new
parametrization, a buy and sell order that were previously both placed at the
price $k$ are now placed at the prices $2k+2$ and $2k+1$, respectively, and
hence still match. We let $\Xc^-_t(2k+2)$ (resp.\ $\Xc^+_t(2k+1)$) denote the
number of buy (resp.\ sell) limit orders in the order book at a given time and
price. We define demand and supply functions
\be\label{ladis}
\la_-:\{3,5,\ldots,2n+1\}\to\R
\quand
\la_+:\{2,4,\ldots,2n\}\to\R
\ee
in such a way that $\la_-(2k+1)$ (resp.\ $\la_+(2k)$) is the total rate at
which buy (resp.\ sell) orders are placed at prices in $\{2k+2,\ldots,2n+2\}$
(resp.\ $\{1,\ldots,2k-1\}$). In particular, $\la_-(2n+1)$ and $\la_+(2)$ are
the rates of buy and sell market orders, respectively.

For any function of the form $f:\{k,k+2,\ldots,m\}\to\R$, we define a discrete
derivative $\di f:\{k+1,k+3,\ldots,m-1\}\to\R$ by
\be
\di f(x):=f(x+1)-f(x-1).
\ee
For sums over sets of the form $\{k,k+2,\ldots,m\}$, we use the shorthand
\be
\sum_k^mg:=\sum_{x\in\{k,k+2,\ldots,m\}}g(x),
\ee
and we define $\sum_k^{k-2}g:=0$. We let
\be
f'(x):=f(x+1)\quand f^\ast(x):=f(x-1)
\ee
denote the function $f$ shifted by one to the left or right, respectively.
It is straightforward to prove the product rule
\be\label{disprodrule}
\di(fg)=f'\di g+g^\ast\di f.
\ee
We also have the following special case of the chain rule:
\be\label{dischainrule}
\di\Big(\frac{1}{f}\Big)=-\frac{\di f}{f'f^\ast}.
\ee

Let
\be\label{diswdom}
w_-:\{2,4,\ldots,2n\}\to\R
\quand
w_+:\{3,5,\ldots,2n+1\}\to\R
\ee
be weight functions satisfying $w_-(2):=0$ and $w_+(2n+1):=0$,
and define a linear function $F$ by
\be\label{disFlin}
F(\Xc):=\sum_{x=4}^{2n}w_-(x)\Xc^-(x)+\sum_{x=3}^{2n-1}w_+(x)\Xc^+(x),
\ee
Then, in analogy with Lemma~\ref{L:GF}, one can check that
\bc\label{disGF}
\dis GF(\Xc)&=&\dis\sum_{M_-(\Xc)+1}^{2n-1}w_+\di\la_+
-w_-\big(M_-(\Xc)\big)\la_+\big(M_-(\Xc)\big)\\[5pt]
&&\dis-\sum_4^{M_+(\Xc)-1}w_-\di\la_-
-w_+\big(M_+(\Xc)\big)\la_-\big(M_+(\Xc)\big).
\ec

In analogy with Theorem~\ref{T:spec}, one can show that if the rates of market
orders $\la_-(2n+1)$ and $\la_+(2)$ are both positive, then, for each
$z\in\{2,4,\ldots,2n\}$, there exist a unique pair of weight functions
$(w^{z,+}_-,w^{z,+}_+)=(w_-,w_+)$ and a unique constant $f_+(z)\in\R$, such
that the linear functional $F^{z,+}=F$ from (\ref{disFlin}) satisfies
\be\label{disfpdef}
GF(\Xc)=1_{\txt\{M_+(\Xc)>z\}}-f_+(z).
\ee
Also, for each $z\in\{3,5,\ldots,2n+1\}$, there exist a unique pair of
weight functions $(w^{z,-}_-,w^{z,-}_+)=(w_-,w_+)$ and constant $f_-(z)$ such
that the linear functional $F^{z,-}=F$ from (\ref{disFlin}) satisfies
\be\label{disfmdef}
GF(\Xc)=1_{\txt\{M_-(\Xc)<z\}}-f_-(z).
\ee

It is possible to derive nice, explicit formulas for these weight
functions. In analogy with (\ref{Gamma}), define
\be\label{disGamma}
\Ga:=\frac{1}{\la_-(2n+1)\la_+(2n)}
-\sum_3^{2n-1}\frac{1}{\la_-}\di\Big(\frac{1}{\la_+}\Big)
=\frac{1}{\la_-(3)\la_+(2)}
+\sum_{4}^{2n}\frac{1}{\la_+}\di\Big(\frac{1}{\la_-}\Big),
\ee
where the equality of both formulas follows from the product rule
(\ref{disprodrule}) applied to the functions $1/\la'_-$ and $1/\la_+$.
Set $I_0:=\{2,4,\ldots,2n\}$ and $I_1:=\{3,5,\ldots,2n+1\}$.
In analogy with (\ref{Vpm}), define
\be\ba{r@{\,}c@{\,}ll}\label{Upm}
\dis u_{-+}(x)&:=&\dis\Ga^{-1}\Big\{
\frac{1}{\la_-(2n+1)\la_+(2n)}
-\sum_{x+1}^{2n-1}\frac{1}{\la_-}\di\Big(\frac{1}{\la_+}\Big)\Big\}
\qquad&\big(x\in I_0\big),\\[5pt]
\dis u_{+-}(x)&:=&\dis\Ga^{-1}\Big\{
\frac{1}{\la_-(3)\la_+(2)}
+\sum_4^{x-1}\frac{1}{\la_+}\di\Big(\frac{1}{\la_-}\Big)
\Big\}\qquad&\big(x\in I_1\big).
\ec
Then, in analogy with Lemma~\ref{L:wint}, one has
\be\ba{r@{\,}c@{\,}ll}\label{diswint}
\dis w^{z,-}_+(x)&=&\dis\la_-(z)\Ga\big[u_{+-}(x\vee z)-1\big]u_{+-}(x\wedge z)
\!&\big(x,z\in I_1\big)\\[5pt]
\dis w^{z,-}_-(x)&=&\dis\la_-(z)\Ga\big(u_{+-}(z)-1_{\{x<z\}}\big)
\big(u_{-+}(x)-1_{\{x<z\}}\big)
\!&\big(x\in I_0,\ z\in I_1\big)\\[5pt]
\dis w^{z,+}_-(x)&=&\dis\la_+(z)\Ga\big[u_{-+}(x\wedge z)-1\big]u_{-+}(x\vee z)
\!&\big(x,z\in I_0\big)\\[5pt]
\dis w^{z,+}_+(x)&=&\dis\la_+(z)\Ga\big(u_{-+}(z)-1_{\{x>z\}}\big)
\big(u_{+-}(x)-1_{\{x>z\}}\big)
\!&\big(x\in I_1,\ z\in I_0\big).
\ec

Moreover, the functions
\be\label{fdis}
f_-:\{3,5,\ldots,2n+1\}\to\R
\quand
f_+:\{2,6,\ldots,2n\}\to\R
\ee
from (\ref{disfpdef}) and (\ref{disfmdef}) satisfy the discrete version of
Luckock's equation, which reads
\be\ba{rr@{\,}c@{\,}ll}\label{disLuckock}
{\rm(i)}&\dis f_-\di\la_+&=&\dis-\la_-\di f_+
&\qquad\mbox{on }\{3,5,\ldots,2n-1\},\\[5pt]
{\rm(ii)}&\dis f_+\di\la_-&=&\dis-\la_+\di f_-
&\qquad\mbox{on }\{4,6,\ldots,2n\},\\[5pt]
{\rm(iii)}&\dis f_+(2)&=&1=f_-(2n+1).
\ec
The solution to this equation can explicitly be written as
\be\ba{rr@{\,}c@{\,}ll}\label{disLucsol}
{\rm(i)}&\dis\Big(\frac{f_+}{\la_+}\Big)(x)
&=&\dis\frac{1}{\la_+(2)}
+\sum_{3}^{x-1}\frac{\kappa}{\la_-}\,\di\Big(\frac{1}{\la_+}\Big)
&\qquad\big(x\in I_0\big),\\[5pt]
{\rm(ii)}&\dis\Big(\frac{f_-}{\la_-}\Big)(x)
&=&\dis\frac{1}{\la_-(2n+1)}
-\sum_{x+1}^{2n}\frac{\kappa}{\la_+}\,\di\Big(\frac{1}{\la_-}\Big)
&\qquad\big(x\in I_1\big),
\ec
where $\kappa$ is given by
\be\label{disLuckappa}
\kappa=\kappa_{\rm L}:=\Ga^{-1}\Big(\frac{1}{\la_-(2n+1)}+\frac{1}{\la_+(2)}\Big),
\ee
and $\Ga>0$ is the constant from (\ref{disGamma}).

\subsection{Suggestions for future work}\label{A:future}

Several open problems concerning Stigler-Luckock models remain. In particular,
these include:
\begin{itemize}
\item[I.] If a Stigler-Luckock model has a competitive window
  $(J^\cri_-,J^\cri_+)$, then show that in the long run, buy orders below
  $J^\cri_-$ and sell orders above $J^\cri_+$ are never matched.
\item[II.] Show that all orders inside the competitive window are eventually
  matched.
\item[III.] Assuming (A3) and (A6), if the solution to Luckock's
  equation satisfies $f_-(I_-)\wedge f_+(I_+)=0$, then show that there
  is an invariant law on $\Si_{\rm ord}$.
\end{itemize}
Under certain technical assumptions, Problems~I and II have been solved in
\cite{KY16}, but Problem~III remains completely open. In view of this, and
also since the methods of \cite{KY16} fail for certain extensions of the model
such as the one discussed in \cite{PS16}, let us investigate how our methods
could possibly be applied to these problems.

Let $(X_k)_{k\geq 0}$ be a Stigler-Luckock model on an interval $I$ and let
$X_k\big|_Z$ be the restriction of $X_k$ to a subinterval $Z=(z_-,z_+)\sub
I$. (Note that this is not what we have called the restricted model on $Z$; in
particular, the latter is a Markov chain, while $X_k\big|_Z$ is not.) To solve
Problem~I, one would need to show that if $Z$ is slightly larger than the
competitive window, then the process $X_k\big|_Z$ is transient, in a suitable
sense, while Problems~II and III could be solved if one could show that if $Z$
is slightly smaller than the competitive window, then the process $X_k\big|_Z$
spends a positive time in the empty state, with some uniform bounds on the
expected number of buy and sell orders in $Z$.

In this context, it is natural to look at linear functions $F$ as in
(\ref{Flin}) such that the weight functions $w_\pm$ are supported on $\ov Z$
and $GF(\Xc)$ depends only on the relative order of $M_\pm(\Xc)$ and $z_\pm$.
It appears that such weight functions exist and form a two-dimensional space.
Using notation as in Theorem~\ref{T:spec}, let us define
\be
\hat w_\pm:=w^{z_-,-}+w^{z_+,+}
\quand
\dot w_\pm:=w^{z_+,-}+w^{z_-,+}.
\ee
Then it appears that there exists a unique constant $c\in\R$ such that
\be
\ov w_\pm:=\hat w_\pm+c\dot w_\pm
\ee
are supported on $\ov Z$. Moreover, it seems that the two-dimensional space 
we just mentioned is spanned by the ``symmetric'' weight functions $(\ov
w_-,\ov w_+)$ and the ``asymmetric'' weight functions $(w^\ast_-,w^\ast_+)$
defined as
\be
w^\ast_-:=-1_{(z_-,z_+]}\quand w^\ast_+:=1_{[z_-,z_+)}.
\ee
Letting $\ov F$ and $F^\ast$ denote the corresponding linear functions,
a natural way to attack Problem~I is to show that if the competitive window $J^\cri$
satisfies $\ov J^\cri\sub Z$, then there exists a function $h(\ov F,F^\ast)$ that
is subharmonic for the generator $G$ in (\ref{Gdef}) and that shows that
$\ov F(\Xc_t)\to\infty$ a.s.\ in such a way that $|F^\ast(\Xc_t)|\ll\ov
F(\Xc_t)$.

Also, a natural way to attack Problems~II and III is to find a ``Lyapunov
style'' function $\Ly$ that depends on $\ov F$, $F^\ast$, and perhaps some other
functions of the process, and that solves an inequality of the form
(\ref{GLya}).

To conclude the paper, we mention a few more open problems, for which we have
nothing more concrete to say.

\begin{itemize}
\item[IV.] Show that the invariant law from Problem~III is unique and the
  long-time limit law started from any initial law.
\item[V.] Investigate existence and uniqueness of solutions to Luckock's
  equation with assumption (A6) replaced by the weaker (A5) plus perhaps some
  conditions involving the function $\Phi$ from (\ref{Phidef}).
\item[VI.] Investigate whether the restricted model on the competitive window
  is null recurrent or transient.\footnote{Theorem~2.1~(1) in the preprint
    \cite{KY16}, which is stated incorrectly, seems to say that it is
    transient. We conjecture that the correct answer is null recurrence,
    however.}
\item[VII.] Prove a limit theorem for the shape of the stationary process near
  the boundary of the competitive window, in the spirit of \cite{FS16}.
\item[VIII.] For the model on the competitive window, investigate the tail of
  the distribution of the time till a limit order is matched.
\end{itemize}

\subsection*{Acknowledgements}

I would like to thank Martin \v{S}m\'id for drawing my attention to Luckock's
paper \cite{Luc03}, and Marco Formentin, Martin Ondrej\'at, and Jan Seidler
for useful discussions. I thank Florian Simatos for bringing the work of Frank
Kelly and Elena Yudovina to my attention, the latter two for answering my
questions about their work, and an unknown referee for a careful reading of
the manuscript.

\end{document}